\documentclass[a4paper,12pt,openright, singlespace]{article}
\usepackage[utf8]{inputenc}
\usepackage{amsmath,amsfonts,amssymb,amsthm}
\usepackage[english]{babel}
\usepackage{graphicx}
\usepackage{epstopdf}
\usepackage{changepage}
\usepackage{subfig}
\usepackage{float}
\usepackage{tikz}
\usepackage[ruled,vlined]{algorithm2e}
\usepackage{algorithmic}
 \usepackage[usenames,dvipsnames]{pstricks}
 \usepackage{epsfig}
 \usepackage{pst-grad} 
 \usepackage{pst-plot} 
 \usepackage[space]{grffile} 
 \usepackage{etoolbox} 
 \makeatletter 
 \patchcmd\Gread@eps{\@inputcheck#1 }{\@inputcheck"#1"\relax}{}{}
 \makeatother
 \usepackage{bm}
 \usepackage{emptypage}

\newcommand{\D}{\displaystyle}

\newcommand{\veps}{\varepsilon}
\newcommand{\dsum}{\D\sum}
\newcommand{\dint}{\D\int}

\newcommand{\mub}{\overline{\mu}}

\newcommand{\uk}{\textbf{u}}

\newcommand{\vk}{\textbf{v}}
\newcommand{\wk}{\textbf{u}}

\newcommand{\wak}{\textbf{w}}

\newcommand{\zk}{\textbf{z}}
\newcommand{\xk}{\textbf{x}}
\newcommand{\unmu}{U_N(\mu)}

\newcommand{\fk}{\textbf{f}}

\newcommand{\intO}[1]{\dint_\Omega #1\; d\Omega}
\newcommand{\intK}[1]{\dsum_{K\in \mathcal{T}_h}\dint_K #1\;d\Omega}
\newcommand{\intk}[1]{\dint_K #1\; d\Omega}

\newcommand{\cD}{\mathcal{D}}
\newcommand{\R}{\mathbb{R}}

\newcommand{\ra}{\rightarrow}

\newcommand{\normlio}[1]{\|#1\|_{0,\infty,\Omega}}

\newcommand{\normlc}[1]{\|#1\|_{0,4,\Omega}}
\newcommand{\normld}[1]{\|#1\|_{0,2,\Omega}}

\newcommand{\normX}[1]{\|#1\|_X}
\newcommand{\normXN}[1]{\|#1\|_N}
\newcommand{\normp}[1]{|||#1|||_N}
\newcommand{\normdual}[1]{\|#1\,\|_{Y'}}
\newcommand{\normtau}[1]{\|#1\|_{\tau_{p}}}

\newcommand{\difuh}{\uk^1_h-\uk^2_h}
\newcommand{\difupih}{\pih\uk^1_h-\pih\uk^2_h}

\newcommand{\difUh}{U^1_h-U^2_h}

\newcommand{\en}{\epsilon_N(\mu)}
\newcommand{\taun}{\tau_N(\mu)}
\newcommand{\dn}{\Delta_N(\mu)}

\newcommand{\nk}{\textbf{n}}
\newcommand{\xip}{\xi^p}
\newcommand{\zetav}{\zeta^\textbf{v}}
\newcommand{\Span}{\mbox{span}}

\newcommand{\Th}{\mathcal{T}_h}

\newcommand{\pih}{\Pi_h^*}
\newcommand{\sigh}{\sigma_h^*}

\newtheorem{lema}{Lemma}
\newtheorem{prop}{Proposition}
\newtheorem{teor}{Theorem}
\newtheorem{corol}{Corollary}
\newtheorem{nota}{Remark}

\newtheorem{hyp}{Hypothesis}

\begin{document}

\title{\LARGE\textbf{On a certified VMS-Smagorinsky Reduced Basis model with LPS pressure stabilization}}

\author{
Tomás Chacón Rebollo\thanks{Departamento de Ecuaciones Diferenciales y Análisis Numérico, 
and Instituto de Matemáticas de la Universidad de Sevilla (IMUS)
, Apdo. de correos 1160, Universidad de Sevilla, 41080 Seville, Spain. chacon@us.es}  
\and Enrique Delgado Ávila\thanks{Departamento de Ecuaciones Diferenciales y Análisis Numérico, Apdo. de correos 1160, Universidad de Sevilla, 41080 Seville, Spain. edelgado1@us.es, macarena@us.es}
\and Macarena Gómez Mármol\footnotemark[2]
}

\date{}

\maketitle

\begin{abstract}
In this work we present the numerical analysis of a Reduced Basis VMS-Smagorinsky model with local projection stabilization (LPS) on the pressure. We construct the reduced velocity space by two different strategies, by considering or not the enrichment of the reduced velocity space with the so-called inner pressure \textit{supremizer}. We present the development of an \textit{a posteriori} error estimator for the snapshot selection through a Greedy algorithm, based on the Brezzi-Rappaz-Raviart (BRR) theory. Moreover, the Empirical Interpolation Method (EIM) is considered for the approximation of the non-linear terms. Finally, we present some numerical tests in which we show an speedup on the computation of the reduced basis problem with the LPS pressure stabilisation, with respect to the method using pressure supremizers.

\end{abstract}

\newpage

\section{Introduction}\label{chap:VMS::sec:introduction}
In this work we present a Reduced Basis (RB) Variational Multi-Scale (VMS) Smago\-rinsky turbulence model with local projection stabilisation (LPS) of the pressure, and small-small VMS modelling of eddy viscosity. We address a turbulence model to consider re\-a\-lis\-tic situations, as turbulent flows frequently appear in actual applications. The Smagorinsky model is a basic Large Eddy Simulation (LES) turbulence model, that provides particularly accurate solutions when the eddy viscosity acts on the resolved small scales (the small-small VMS setting, cf. \cite{LES-VMS}). \black{A way to build ROMs for turbulence models is to derive them from LES models, in particular the Smagorinsky model. This approach has been followed in \cite{PaperSmago} for the Smagorinsky model and in \cite{Delgado2020} for VMS-Smagorinsky model. An alternative consists in initially building a ROM for Navier-Stokes equations, and then model the unresolved scales of the ROM (either by eddy diffusion or other techniques) to build the turbulence model (see  e.g. \cite{iliescu,VMS-Iliescu}).}

Stabilised methods for finite element discretisations provide stable discretisations of the pressure that allow to circumvent the Brezzi-Babuska inf-sup condition. This occurs with both residual-based methods and penalty methods. As residual-based stabilisation methods let us mention the Galerkin-Least Squares (GALS) (cf. \cite{GALS}) and modifications of this method as the Streamline Upwind Galerkin (SUPG) method (cf. \cite{SUPG, Num_approx_book}), or the Adjoint-Stabilised method (cf. \cite{Ad-Stab}). The stabilising effect of those methods is due to the dissipative effect of the terms added to the Galerkin projection. In particular, the pressure stabilisation is granted thanks to a weighted pressure Laplacian term. The penalty methods are simplified versions that just keep the dissipative terms, although loosing full consistency. We may mention the penalty term-by-term stabilised method (cf. \cite{TbT_Tomas98}), that is an extension of the pure penalty method introduced in \cite{Stab_penalty}. In the latter method, the penalty term is just a scaled pressure Laplacian. Further, the Local Projection Stabilisation (LPS) treatment of the pressure is a penalty method that just stabilises the high-frequency components of the pressure gradient that cannot be represented in the discrete velocity space. It allows to use non-stable pairs of velocity-pressure finite elements, such as $\mathbb{P}_k-\mathbb{P}_k$, $k \ge 2$, while keeping optimal order of accuracy in both velocity and pressure (cf. \cite{TbT_Tomas98, Tomas_TbTStab}).

There exist several well-established techniques to stabilise the dis\-cre\-ti\-sa\-tion of the reduced pressure in ROM of in\-com\-pre\-ssi\-ble flows. Let us mention to enrich the reduced velocity space with reduced velocity divergence - reduced pressure supremizers (cf. \cite{Ballarin, Stokes3}), or solving a Poisson equation for the reduced pressure --by testing the momentum conservation equation with gradients of reduced pressures--({cf. \cite{caiazzo2014numerical, Patera}), or using SUPG stabilisation (cf. \cite{Ali2020, Novo2022}).  
These techniques provide efficient approximations of the reduced pressure, although are somewhat involved, requiring either increasing the dimension of the reduced velocity space, building an equation for the pressure with a number of terms, or increasing the complexity of the discrete problem. The LPS stabilisation provides a simplified alternative to these techniques.  The POD-ROM solution of  Navier-Stokes equations with LPS stabilisation of the pressure has recently been addressed in \cite{Novo2021}, which includes an a-priori error analysis of the method. \black{Other works with a priori error analysis for VMS ROM based on POD are reported in \cite{Eroglu2017, Rubino2022}. Moreover, further works on stabilised ROM for convection-dominated problems can be found in \cite{Giere2015, Stab_RBM} where a SUPG stabilization is considered for both full order model (FOM) and ROM, with an offline/online structure for the ROM, whereas in \cite{Bergmann2009} different stabilisation procedures, as SUPG and VMS, are performed only on the POD-ROM using the residual of the Navier-Stokes operator, evaluated on the discarded POD modes.}

%

In this work, we address the \textit{a posteriori} error analysis - based reduced order modelling of incompressible flow equations with LPS stabilisation of the pressure. This allows to construct \lq\lq certified" solvers, with error level below targeted values. We actually afford the reduced basis solution of the steady Smagorinsky turbulence model for the sake of adressing realistic flows.  We include the small-small VMS formulation of the eddy diffusion which, ac\-tua\-lly, is a three-level LES turbulence model ({cf. \cite{BernardiVMS}) that provides good accuracy in the approximation of first and second order turbulence statistics ({cf. \cite{VMS_numAnal}). We prove the stability of the reduced problem with LPS stabilisation of the pressure for velocity-pressure finite elements, whenever the pressure is piecewise affine. We construct an \textit{a posteriori} error bound estimator u\-sing the Brezzi-Raviart-Rappaz (BRR) theory of approximation of regular branches of solutions of non-linear equations (cf. \cite{BRR}). This is the key to build a certified RB solver. We use the Empirical Interpolation Method (EIM) to build reduced approximations of the non-linear coefficients appearing in model: the eddy diffusion and the statibilised coefficient appearing in the LPS formulation. We present some numerical results \black{where} we compare the velocity supremizers approach versus the LPS pressure discretization. We observe an improved speedup for the LPS discretization with similar error levels.

The structure of this work is as follows. In section \ref{chap:VMS::sec:LPS-VMS}, we present the VMS-Smagorinsky model with local projection stabilisation for the pressure. In subsection \ref{chap:VMS::sec:HighFid}, we present the finite element (FE) model, while the reduced basis method is considered in subsection \ref{chap:VMS::sec:RBP}. The numerical analysis for the well-posedness is presented in section \ref{chap:VMS::sec:well-posedness}, and the sub-sequent development of the \textit{a posteriori} error bound estimator in section \ref{chap:VMS::sec:Post_Err}. Then, in section \ref{chap:VMS::sec:EIM}, we present the approximation of the eddy-viscosity and the pressure stabilisation constant by the EIM.  Finally, we present the numerical tests in section \ref{chap:VMS::sec:num_results}, in which we highlight the speedup rates for the computation of solution with the present RB model presented in this work.

\section{Pressure local projection stabilised VMS-Smagorinsky model}\label{chap:VMS::sec:LPS-VMS}
In this section we describe a local projection-stabilised VMS-Smagorinsky model. In the turbulence LES (Large Eddy Simulation) high fidelity model that we consider, the eddy-viscosity only acts on the small resolved scales. This leads to a less diffusive model than the Smagorinsky one presented in \cite{PaperSmago}, where the eddy viscosity acts on both large and small resolved scales. In \cite{VMS_numAnal}, it has been shown that this model has a good precision for the approximation of the mean flow and the second-order statistics of the turbulent flow. Moreover, we use a LPS stabilisation treatment for the pressure, that allows us to use non-stable pairs of finite elements.

\subsection{Finite element problem}\label{chap:VMS::sec:HighFid}
Let $\Omega$ be a bounded domain of $\R^d\, (d=2,3)$, with Lipschitz-continuous boundary $\Gamma$. We assume that $\Gamma$ is split into $\Gamma=\Gamma_D \cup \Gamma_N$ where $\Gamma_D$ and $\Gamma_N$ are two connected measurable sets of positive $(d-1)$-dimensional measure, with disjoint interiors. We intend to impose Dirichlet and Neumann boundary conditions on $\Gamma_D$ and $\Gamma_N$, respectively.

We present a parametric steady Smagorinsky model, where we consider the Reynolds number as a parameter, denoted by $\mu\in\cD$, where $\cD$ is a compact sub-set of $\R$. Although the Smagorinsky model is intrinsically discrete, we present it in a continuous form in order to clarify its relationship with the Navier-Stokes equations: We search for a velocity field $\wk \,: \Omega \mapsto \R^d$ and a pressure function $p\,:\Omega \mapsto \R$ such that
\begin{equation}\label{NS}\left\{\begin{array}{ll}
\wk\cdot\nabla\wk+\nabla p-\nabla\cdot\left(\dfrac{1}{\mu}\nabla\wk+\nu_T(\wk')\nabla\wk'\right)=\fk&\mbox{ in }\Omega,\vspace{0.2cm}\\

\nabla\cdot\wk={\bf 0}&\mbox{ in } \Omega,\vspace{0.1cm}\\
\wk=0&\mbox{ on }\Gamma_D,\vspace{0.1cm}\\
\nk\cdot\left(\dfrac{1}{\mu}\nabla\wk+\nu_T(\wk')\nabla\wk'\right)=0&\mbox{ on }\Gamma_N,\vspace{0.2cm}\\
\end{array}\right.
\end{equation}
where $\fk$ is the kinetic momentum source, $\nu_T(\wk')$ is the eddy viscosity defined as
\begin{equation} \label{eq:eddivis}
\nu_T(\wk') = C_S^2\D\sum_{K\in\mathcal{T}_h}h_K^2\big|\nabla\wk_{|_K}'\big|\chi_K, 
\end{equation} 
where $\big|\cdot\big|$ denotes the Frobenius norm in $\R^{d\times d}$, $C_S$ is the Smagorinsky cons\-tant and $\wk'$ stands for the small-scales part of the resolved velocity field $\wk$, $\wk'$ is assumed to lie in the inertial spectrum of the full velocity field. The eddy viscosity setting in \eqref{NS}-\eqref{eq:eddivis} corresponds to the small-small VMS setting introduced in  \cite{LES-VMS}. Other possibilities would correspond to replacing $\wk'$ by $\wk$ in either \eqref{NS} or \eqref{eq:eddivis}. The eddy viscosity in \eqref{NS} only acts on the small scales of the resolved part of the velocity field. Let $\{\Th\}_{h>0}$ is uniformly regular mesh in the sense of Ciarlet \cite{Ciarlet1978}.

To state the full order discretization that we consider for problem \eqref{NS}, let us introduce the spaces
$$
Y=\{ \vk \in H^1(\Omega)^d,\,\, \mbox{s.t. } \vk_{|_{\Gamma_D}}={\bf 0}\,\},\quad M=\{q \in L^2(\Omega),\,\, \mbox{s.t. }\, \int_\Omega q =0\,\}.
$$
We assume $\fk \in Y'$. 

Given an integer $n\ge 1$ and an element $K\in\Th$, we denote by $\mathbb{P}_n(K)$ the Finite Element space given by Lagrange polynomials of degree less than, or equal to, $n$ defined on $K$. We assume that $\Omega$ is polygonal when $d=2$ or polyhedric when $d=3$, and consider a triangulation $\Th$ such that $\Gamma_D$ and $\Gamma_N$ are unions of full simplices of the boundaries of elements of $\Th$. Let us define the finite element space
\begin{equation}\label{FE_Space_Vn}
V^{(n)} = \{v_h \in C^0(\overline{\Omega}) \text{ such that } v_{h|_K}\in\mathbb{P}_n(K), \forall K\in\Th \} .
\end{equation}
Let $Y_h = Y \cap [V^{(n)}]^d $ and $M_h \subset M \cap V^{(m)}$ be two finite element subspaces of $Y$ and $M$, respectively.  We introduce the LPS-VMS discretisation of the Smagorinsky model \eqref{NS},  with stabilisation on the pressure, 
\begin{equation}
\label{chap:VMS::pb:FV}\left\{\begin{array}{l}
\mbox{Find } (\uk_h,p_h)=(\uk_h(\mu),p_h(\mu))\in Y_h\times M_h\mbox{ such that}\vspace{0.3cm}\\
\begin{array}{ll}
a(\uk_h,\vk_h;\mu)+b(\vk_h,p_h;\mu)+ a_S'(\wk_h;\wk_h,\vk_h;\mu) \\
+c(\uk_h,\uk_h,\vk_h;\mu)=\left<\fk,\vk\right>&\quad\forall\vk_h\in Y_h,\vspace{0.1cm}\\
b(\uk_h,q_h;\mu)+s_{pres}(p_h,q_h;\mu)=0&\quad\forall q_h\in M_h,\end{array}\end{array}\right.
\end{equation} 
where the bilinear forms $a(\cdot,\cdot;\mu)$ and $b(\cdot,\cdot;\mu)$ are defined as 
\begin{equation}\label{chap:Smago::eq:ab_form}
a(\uk,\vk;\mu)=\frac{1}{\mu}\int_\Omega\nabla\uk:\nabla\vk\,d\Omega,\qquad b(\vk,q;\mu)=-\int_\Omega(\nabla\cdot\vk)q\,d\Omega;
\end{equation}
while the trilinear form $c(\cdot,\cdot,\cdot;\mu)$ is defined as
\begin{equation}\label{chap:Smago::eq:conv_form}
c(\zk,\uk,\vk;\mu)=\frac{1}{2}\, \left [\int_\Omega(\zk\cdot\nabla\uk)\vk\,d\Omega- \int_\Omega(\zk\cdot\nabla\vk)\uk\,d\Omega\, \right ].
\end{equation}
Moreover, the non-linear form  $a'_S(\cdot;\cdot,\cdot;\mu)$, is a multi-scale Smagorinksy mo\-de\-lling for the eddy viscosity term, and it is given by
\begin{equation}\label{chap:VMS::eq:Small-Small}
a_S'(\zk_h;\uk_h,\vk_h;\mu)=\intO{\nu_T(\pih\zk_h)\nabla(\pih\uk_h):\nabla(\pih\vk_h)},
\end{equation}
where  $\Pi_h^*=Id-\Pi_h$, $\Pi_h$ being a uniformly stable (in $L^2$ and $H^1$-norms) interpolation operator on $Y \cap [V_h^{(n-1)}]^d$. 
This interpolation operator must satisfy optimal error estimates (cf. \cite{desinv} Chap. IX) and preserve the homogeneous Dirichlet boundary conditions when restricted to $Y_h$. In this way the small-scale component of the velocity field that we consider is defined by $\wk_h'=\pih\wk_h$. The eddy viscosity only acts on $\pih\wk_h$. \black{For the numerical experiments, we consider a standard nodal Lagrange operator for its stability and accuracy properties, as well as its simplicity to be computed. (\textit{cf.} \cite{VMS_numAnal}).}

In VMS terminology, the \textit{Large-Small} setting would corresponds to modeling the turbulent viscosity by a function of the whole resolved velocity, taking
\begin{equation}\label{chap:VMS::eq:Small-All}
a_S'(\zk_h;\uk_h,\uk_h)=\intO{\nu_T(\zk_h)\nabla(\pih\uk_h):\nabla(\pih\vk_h)}.
\end{equation} 

Further, the term $s_{pres}(\cdot,\cdot;\mu)$ is the projection-stabilisation term for the pressure, defined as
\begin{equation}\label{chap:VMS::eq:s_pres}
s_{pres}(p_h,q_h;\mu)=\D\sum_{K\in\Th}\tau_{p,K}(\mu)(\sigh(\nabla p_h),\sigh(\nabla q_h))_K.
\end{equation}

Here, $\sigh=Id-\sigma_h$, with $\sigma_h$ some locally stable $L^2$ projection or interpolator operator from $L^2(\Omega)$ to $[V_h^{(n)}]^d$. 
Moreover, $\tau_{p,K}(\mu)$ is the stabilisation coefficient  in (\ref{chap:VMS::eq:s_pres}), that must verify the following hypotesis:
\begin{hyp}\label{chap:VMS::hyp:tau_acot}
There exists two positive constants $\alpha_1,\alpha_2$, independent of $h$, such that
\begin{equation}\label{chap:VMS::eq:estab_coef_acot}
\alpha_1 h_K^2\le\tau_{p,K}(\mu)\le \alpha_2 h_K^2, \quad \forall \mu\in\cD, \forall K\in\Th, \forall h>0.
\end{equation}
\end{hyp}

In the following, we will use the stabilisation pressure coefficient $\tau_{p,K}(\mu)$ proposed by Codina in \cite{Codina_stab}, and used in \cite{VMS_numAnal, phdSamu}:
\begin{equation}\label{chap:VMS::eq:stab_coeff}
\tau_{p,K}(\mu)=\left[c_1 \dfrac{1/\mu+\overline{\nu}_{T_{|K}}(\mu)}{h_K^2}+c_2\dfrac{U_K(\mu)}{h_K} \right]^{-1},
\end{equation}
where $\overline{\nu}_{T_{|K}}$ is some local eddy viscosity, $U_K$ is a local velocity and $c_1, c_2$ some positive experimental constants. Taking $\tau_{p,K}(\mu)$ by this way, we are ensuring (\ref{chap:VMS::eq:estab_coef_acot}).

By denoting $X_h=Y_h\times M_h$, we rewrite problem ($\ref{chap:VMS::pb:FV}$), with compact notation, as follows
\begin{equation}\label{chap:VMS::pb:FVX}\left\{\begin{array}{l}
\mbox{Find }U_h(\mu)=(\uk_h,p_h)\in X_h \mbox{ such that}\vspace{0.2cm}\\
A(U_h(\mu),V_h;\mu)=F(V_h;\mu)\qquad\forall V_h\in X_h,\end{array}\right.
\end{equation}
where
\begin{equation}\label{chap:VMS::aff_decomp}
\begin{array}{ll}
A(U_h(\mu),V_h;\mu)&=\dfrac{1}{\mu}A_{Diff}(U_h,V_h)+A_{Div}(U_h,V_h)+A_{Conv}(U_h;V_h)\vspace{0.1cm}\\
&+A_{Sma}(U_h;V_h)+A_{Pres}(U_h,V_h;\mu),
\end{array}
\end{equation}
where, denoting $U=(\uk,p^u)$, $V=(\vk,p^v)$, 
\[
\begin{array}{rll}
A_{Diff}(U,V)&\!\!\!\!=\intO{\nabla\uk:\nabla\vk},\vspace{0.2cm}\\                                                
A_{Div}(U,V)&\!\!\!\!=\intO{(\nabla\cdot\uk )p^v}-\intO{(\nabla\cdot\vk)p^u}\vspace{0.1cm}\\
A_{Conv}(U;V)&\!\!\!\!=\displaystyle\frac{1}{2}\,\intO{\left ( (\uk\cdot\nabla\uk)\vk- (\uk\cdot\nabla\vk)\uk \,\right )},\vspace{0.2cm}\\
A_{Sma}(U;V)&\!\!\!\!=\intO{\nu_T(\pih\wk)\,\nabla(\pih\wk):\nabla(\pih\vk)},\vspace{0.2cm}\\
A_{Pres}(U,V;\mu)&\!\!\!\!=\D\sum_{K\in\Th}\intk{\tau_{p,K}(\mu)\,\sigh(\nabla p^u)\cdot\sigh(\nabla p^v)},\vspace{0.2cm}\\
F(V_h;\mu)&\!\!\!\!= \langle \fk, \vk_h \rangle.
\end{array}
\] 

The numerical analysis (existence of solutions, stability and error analysis) of the full-order projection-based VMS-Smagorinsky model for problem \eqref{chap:VMS::pb:FV} (or equivalently \eqref{chap:VMS::pb:FVX}) with wall-law boundary conditions can be found in \cite{VMS_numAnal}. 
\subsection{Reduced basis problem}\label{chap:VMS::sec:RBP}
In this section we introduce the Reduced Basis (RB) model. It is the discretisation (\ref{chap:VMS::pb:FVX}) with LPS stabilisation of the pressure built on the reduced spaces.  The RB problem reads 
\begin{equation}\label{chap:VMS::pb:RBP}\left\{\begin{array}{l}
\mbox{Find } (\uk_N,p_N)\in Y_N\times M_N\mbox{ such that}\vspace{0.3cm}\\

\begin{array}{ll}
a(\uk_N,\vk_N;\mu)+b(\vk_N,p_N;\mu)+ a_S'(\wk_N;\wk_N,\vk_N;\mu) \\
+c(\uk_N,\uk_N,\vk_N;\mu)=\langle \fk, \vk_N\rangle &\quad\forall\vk_N\in Y_N,\vspace{0.1cm}\\
b(\uk_N,q_N;\mu) \black{ + s_{pres}(p_N,q_N;\mu)}=0&\quad\forall q_N\in M_N.\end{array}\end{array}\right.
\end{equation} 

Here, %
we denote by $Y_N$ the reduced velocity space, and by  $M_N$ the reduced pressure space. The computation of the reduced spaces is done through the Greedy algorithm \cite{Greedy70s}. The Greedy algorithm summarizes: 
\begin{enumerate}
\item Set $\mu^1$, compute $(\uk_h(\mu^1), p_h(\mu^1))$ solving problem (\ref{chap:VMS::pb:FV}), and define the reduced space $Y_1=\Span\{\uk_h(\mu^1)\}$, $M_1=\Span\{p_h(\mu^1)\}$. Orthonormalize the spaces $Y_1, M_1$.
\item For $k\ge2$, compute $\Delta_{k-1}(\mu), \forall\mu\in\cD_{train}$ and set $\mu^{k}=\arg\D\max_{\mu\in\cD_{train}}\Delta_{k-1}(\mu).$ Here, we denote by $\Delta_k$ the \textit{a posteriori} error estimator.
\item Compute $(\uk_h(\mu^{k}), p_h(\mu^k))$, define the reduced spaces $Y_k, M_K$ by adding the new computed snapshot, and orhtonormalize them.
\item Stop if $\D\max_{\mu\in\cD_{train}}\Delta_k(\mu)<\varepsilon_{RB}$, \black{with  $\veps_{RB}$ a prescribed tolerance for the Greedy algorithm} . If not, back to 2.
\end{enumerate} 

This yields reduced velocity and pressure spaces corresponding to the solution of problem \eqref{chap:VMS::pb:RBP} for selected values of the Reynolds number:
\begin{equation}\label{chap:VMS::eq:vel_space2}
Y_N=\mbox{span}\{\zetav_i:=\uk(\mu^i),\; i=1,\dots,N\};
\end{equation}
\begin{equation}
M_N=\mbox{span}\{\xip_i:=p(\mu^i),\; i=1,\dots,N\},
\end{equation}
where the $\mu^i$ are selected by the Greedy algorithm. The stability of this discretization is proved in the next section.

To test the numerical performances of this discretisation, we compare it with the strategy consisting in  enriching the reduced velocity space with the supremizers of the reduced velocity divergence - reduced pressure duality.  The reduced velocity space in this case is given by
\begin{equation}\label{chap:VMS::eq:vel_space}
Y_N=\mbox{span}\{\zetav_{2i-1}:=\uk(\mu^i),\zetav_{2i}:=T_p^\mu(\xip_i)),\; i=1,\dots,N\},
\end{equation}
where here we are denoting $T_p^\mu:M_h\ra Y_h$ the reduced pressure \textit{supremizer}, defined by 
\begin{equation}\label{supremizer_def}
\left(\nabla T_p^\mu q_h,\nabla \vk_h\right)_\Omega=b(q_h,\vk_h;\mu)\quad\forall\vk_h\in Y_h.
\end{equation}

This strategy leads to an offline/online stabilisation where the dimension of the reduced space $X_N=Y_N \times M_N$ in problem \eqref{chap:VMS::pb:RBP} is $3N$, while using the LPS stabilisation the dimension of $X_N$ is $2N$.

When $Y_N$ is defined by  \eqref{chap:VMS::eq:vel_space}, i. e., when it includes the reduced pressure supremizers, then the discrete inf-sup condition between spaces $Y_N$ and $M_N$ holds  (cf. \cite{Ballarin}). Consequently, for small enough data $\fk$, the reduced problem \eqref{chap:VMS::eq:vel_space2} admits a solution by a finite-dimensional compactness argument, similar to the one used in \cite{VMS_numAnal} to prove the existence of solutions of high fidelity problem \eqref{chap:VMS::pb:FV}. 

When the LPS approach is followed, this argument also holds when the \black{pressures} are piecewise affine. Indeed, in this case we may bound a weaker norm of the pressure,
\begin{equation}\label{infsupred}
\normp{q_N}=\sup_{\vk_N\in Y_N}\dfrac{(q_N,\nabla\cdot\vk_N)_\Omega}{\normld{\nabla{\vk_N}}}+\normtau{\sigh(\nabla q_N)}, \quad \forall q_N\in M_N,
\end{equation}
where
$$
\normtau{r}=\left (\D\sum_{K\in\Th}\tilde{\tau}_{p,K}\|r\|_{0,2,K}^2 \right)^{1/2},\quad \mbox{for }\,\, r \in L^2(\Omega),
$$
$\mbox{where  }\,\, \tilde{\tau}_{p,K}=\displaystyle\inf_{\mu \in \cD}\tau_{p,K}(\mu)$, and we denote by $\|\cdot\|_{l,p,S}$ the $W^{l,p}(S)$ norm on some open set $S$.
It holds
\begin{prop} \label{propnormp}
Assume that $M_h$ is formed by piecewise affine functions, and that the space $Y_N$ defined by \eqref{chap:VMS::eq:vel_space2} contains some function $\wak_N$ such that $ \int_\Omega \wak_N \ne 0$. Then $\normp{\cdot}$ is a norm on $M_N$. Moreover, there exists a constant $C>0$ such that
\begin{equation} \label{boundnormp}
\normp{q_N}\le C\, \|q_N\|_{0,2,\Omega}, \quad \forall q_N\in M_N.
\end{equation}
\end{prop}
\textbf{Proof.} It is straightforward to prove that $\normp{q}\ge 0$ for any $q\in L^2(\Omega)$ and that $\normp{\cdot}$ satisfies the \black{triangle} inequality. Assume now that $\normp{q_N}=0$ for some $q_N \in M_N$. Then $\sigh(\nabla q_N)=0$ and consequently $\nabla q_N=\sigma_h(\nabla q_N) \in [V_h^{n}]^d$. Then $\nabla q_N$ is a continuous function and, as it is a constant on each element $K\in {\cal T}_h$, then  $\nabla q_N$ is constant on $\Omega$. Moreover, as
$\displaystyle
\sup_{\vk_N\in Y_N}\dfrac{(q_N,\nabla\cdot\vk_N)_\Omega}{\normld{\nabla{\vk_N}}}=0,
$
then $(q_N,\nabla\cdot\vk_N)_\Omega=0$ for any $\vk_N \in Y_N$. Then
$$
\nabla q_N \, \int_\Omega \wak_N = - (q_N,\nabla\cdot \wak_N)_\Omega=0.
$$
Consequently, $\nabla q_N=0$ and as $\displaystyle \int_\Omega q_N=0$ because $q_N \in M_h$, it follows $q_N=0$.

Finally, estimate \eqref{boundnormp} readily follows from the inverse finite element estimate (cf. \cite{desinv})
$$
\|\nabla q_N\|_{0,2,K}\le C \, h_K^{-2}\,\|q_N\|_{0,2,K}.
$$
\qed

Then it follows
\begin{corol} Under the hypotheses of Proposition \ref{propnormp}, the reduced problem  \eqref{chap:VMS::pb:RBP} admits at least a solution, that satisfies the estimates
$$
\|\nabla \uk_N\|_{0,2,\Omega} +\normp{p_N} \le \Phi(\|\fk\|_{Y'};\mu).
$$
for some positive continuous function $\Phi(\eta;\mu)$ increasing to $+\infty$ as  $\eta \mapsto +\infty $ or \black{$\mu \mapsto +\infty$}.
\end{corol}
\textbf{Sketch of the proof.} The proof of existence follows from a finite-dimensional compactness argument using Brouwer's fixed point theorem, similar to the one used in \cite{VMS_numAnal} to prove the existence of solutions of the high fidelity problem \eqref{chap:VMS::pb:FV} with wall-law boundary conditions through a linearisation procedure. The essential point is the bounds for velocity and pressure that on one hand  imply the uniqueness of solutions of the linearised problem and on another hand guarantee the hypotheses of Brouwer's fixed point theorem. Indeed, setting $\vk_N=\uk_N$ and $q_N=-p_N$, and summing the two equations appearing in \eqref{chap:VMS::pb:RBP}, it follows
\begin{equation}\label{estunopN}
\|\nabla \uk_N\|_{0,2,\Omega} +\normtau{\sigh(\nabla p_N)} \le \Phi_1(\|\fk\|_{Y'};\mu)
\end{equation}
for some positive continuous function $\Phi_1(\eta;\mu)$ increasing to $+\infty$ as $\eta \mapsto +\infty$ or $\mu \mapsto 0$. Further, setting $q_N=0$ in \eqref{chap:VMS::pb:RBP}, it holds
$$
(\nabla\cdot \vk_N, p_N)_\Omega=a(\wk_N,\wk_N;\mu)+ a_S'(\wk_N;\wk_N,\vk_N;\mu) + c(\wk_N,\wk_N,\vk_N;\mu)-\langle \fk,\vk_N\rangle 
$$ 
for all $\vk_N\in Y_N$. Consequently,
$$
(\nabla\cdot \vk_N, p_N)_\Omega \le \Phi_2(\|\fk\|_{Y'};\mu)\, \|\nabla \vk_N\|_{0,2,\Omega} 
$$
for some function $\Phi_2(\cdot;\cdot)$ with the same properties as $\Phi_1$. Then,  
\begin{equation}\label{estdospN}
\sup_{\vk_N\in Y_N}\dfrac{(\nabla\cdot\vk_N,p_N)_\Omega}{\normld{\nabla{\vk_N}}} \le \Phi_2( \|\fk\|_{Y'};\mu).
\end{equation}
Combining \eqref{estunopN} and \eqref{estdospN} it follows
$$
\|\nabla \uk_N\|_{0,2,\Omega} +\normp{p_N} \le \Phi_1(\|\fk\|_{Y'};\mu)+ \Phi_2( \|\fk\|_{Y'};\mu).
$$
The same argument allows to prove the bounds for velocity and pressure in the linearised problem and then to apply Brouwer's fixed point theorem. We do not describe the proof in detail for brevity.
\qed
\section{Well-posedness analysis}\label{chap:VMS::sec:well-posedness}
In this section we study the well-posedness of the full order \eqref{chap:VMS::pb:FV} and of the reduced problem \black{(\ref{chap:VMS::pb:RBP})}, as well as of the linearised problem. This will be the basis for the derivation of the a posteriori error estimation in the next section.



We shall denote by $C>0$ along the paper the constants that may vary from line to line but which are independent of $\nu$ and $h$. We need the following basic results: 

\begin{lema}\label{estformconv} It holds 
\begin{equation}\label{esformcon}
|c(\zk,\uk,\vk;\mu)| \le C_{S4}^2\, \|\zk\|_{0,2,\Omega}\, \|\uk\|_{0,2,\Omega}\, \|\vk\|_{0,2,\Omega}\,\forall \, \zk,\,\uk,\,\vk \in Y,
\end{equation}
where $C_{S4}$ is the norm of the injection of $Y$ into $L^4(\Omega)^d$.
\end{lema}
This is a standard result, we omit its proof \black{(see e.g. \cite{TomasSmago} Lemma 6.3)}.
\begin{lema}\label{chap:VMS::prop:as_cont}
For any $\uk_h, \vk_h, \zk_h \in Y_h$ and for any $\mu\in\cD$, the non-linear form $a_S'(\cdot;\cdot,\cdot;\mu)$ defined in \eqref{chap:VMS::eq:Small-Small}  satisfies
\begin{equation}\label{chap:VMS::eq:as_cont}
|a_S'(\zk_h;\uk_h,\vk_h;\mu)|\le C_f \,C_S^2h^{2-d/2}\normld{\nabla\zk_h}\normld{\nabla\uk_h}\normld{\nabla\vk_h},
\end{equation}
for some constant $C_f >0$.
\end{lema}
\textbf{Proof.} Using the local inverse finite element estimates (cf. \cite{desinv})
$$
\|\nabla (\pih \zk_h)\|_{0,\infty,K} \le C \, h_K^{-d/2} \,\|\nabla (\pih \zk_h)\|_{0,2,K}
$$
for some constant $C>0$, it holds $|\nu_T(\pih \zk_h) ({\bf x})|\le C\,C_S^2\, h_K^{2-d/2}\,\|\nabla (\pih \zk_h)\|_{0,2,K}$ for all ${\bf x} \in K$, for any $K \in {\cal T}_h$ . Then,
\begin{eqnarray*}
&&|a_S'(\zk_h;\uk_h,\vk_h;\mu)|  \\ 
&&\le C\,C_S^2 \, h^{2-d/2}\,\|\nabla (\pih \zk_h)\|_{0,2,\Omega}\, \sum_{K \in {\cal T}_h}\|\nabla (\pih \uk_h)\|_{0,2,K}\, \|\nabla (\pih \vk_h)\|_{0,2,K}\\
&&\le C\,C_S^2 \, h^{2-d/2}\|\nabla (\pih \zk_h)\|_{0,2,\Omega} \,\|\nabla (\pih \uk_h)\|_{0,2,\Omega}\, \|\nabla (\pih \vk_h)\|_{0,2,\Omega}.
\end{eqnarray*}
Using the stability of operator $\Pi_h$ in $H^1(\Omega)^d$ norm, estimate \eqref{chap:VMS::eq:as_cont} follows.

\qed

\begin{lema}\label{chap:VMS::lema:eq_norm_tau}
Assume Hypothesis \ref{chap:VMS::hyp:tau_acot} holds, and let $q_h\in V_h^l(\Omega)$. Then,
\begin{equation}\label{chap:VMS::eq:des_norm_tau}
\normtau{\sigh(\nabla q_h)}\le C_{\tau} \normld{q_h}
\end{equation}
for some constant $C_{\tau} >0$.
\end{lema}
\textbf{Proof.} Taking into account Hypothesis \ref{chap:VMS::hyp:tau_acot}, and the local stability of operator $\sigma$, it holds
\begin{equation}\label{chap:VMS::eq:des_nablaq}
\tau_K \|{\sigh(\nabla q_h)}\|_{0,2,K}^2\le  C h_K^2 \|{\nabla q_h}\|_{0,2,K}^2
\end{equation}
Considering the local inverse inequalities for finite element functions (cf. \cite{desinv}) it follows
\begin{equation}
\|{\nabla q_h}\|_{0,2,K} \le C\, h_K^{-1}\,\|q_h\|_{0,2,K},\,\, \forall K \in {\cal T}_h.
\end{equation}
Consequently,
\[
\normtau{\sigh(\nabla q_h)}\le C \,\normld{q_h}.
\]
\qed

%

We shall use the following inf-sup generalised condition (cf. \cite{Tomas_TbTStab})
\begin{prop}\label{chap:VMS::prop:infsup_b}
Assume that the family of triangulations $\{\Th\}_{h>0}$ is uniformly regular. Then it holds
\begin{equation}\label{eq:infsup_b}
\normld{q_h}\le \alpha\left(\sup_{\vk_h\in Y_h}\dfrac{(q_h,\nabla\cdot \vk_h)_\Omega}{\normld{\nabla\vk_h}}+\normtau{\sigh(\nabla q_h)}\right) \quad\forall q_h\in M_h,
\end{equation}
with $\alpha>0$, independent of $h$.
\end{prop}
We will assume in the sequel for simplicity of derivations that the grids are uniformly regular. The previous result can be extended to regular rather than uniformly regular grids, with a more involved analysis. For further details, see \cite{Tomas_TbTStab}.

The following \textit{a priori} estimates for the full order velocity and pressure hold:
\begin{teor}\label{chap:VMS::prop:apriori_estimates}
Let $(\uk_h,p_h)\in X_h$ a solution of problem \eqref{chap:VMS::pb:FV}. Let $\fk\in Y_h'$. Then, the following \textit{a priori} estimates are verified:
\begin{equation}\label{chap:VMS::eq:estimate1}
\normld{\nabla\wk_h}\le \mu \,C_P\,\normdual{\fk},
\end{equation}
\begin{equation}\label{chap:VMS::eq:estimate2}
\normtau{\sigh(\nabla p_h)}\le \sqrt{\mu} \,C_P\,\normdual{\fk},
\end{equation}
\begin{equation}\label{chap:VMS::eq:estimate3}
\|p_h\|_{0,2,\Omega}\le \black{\alpha}\,\left ( ((1+\sqrt{\mu})\,C_P+1)\,\normdual{\fk}+\frac{C}{\mu\, C_P} \, \normdual{\fk}^2 \right ),
\end{equation}
where $C_P>0$ is the norm of the Poincar\'e injection of $Y$ into $L^2(\Omega)$ and $C>0$ is a constant independent of $h$, $\mu$ and the data $\fk$.
\end{teor} 
\textbf{Proof.} Estimates \eqref{chap:VMS::eq:estimate1} and \eqref{chap:VMS::eq:estimate2} are readily obtained taking $\vk_h=\wk_h$ and $q_h=p_h$ in \eqref{chap:VMS::pb:FV} and summing the two identities appearing therein.
%

To prove \eqref{chap:VMS::eq:estimate3}, we take $\vk_h=\wk_h$ and $q_h=0$ in \eqref{chap:VMS::pb:FV}. We have
\[
a(\wk_h, \vk_h;\mu) + a_s'(\wk_h;\wk_h,\vk_h;\mu) + c(\wk_h,\wk_h,\vk_h;\mu) - (p_h, \nabla\cdot\vk_h) = \left<\fk, \vk_h\right>.
\]
Thus, from Cauchy-Schwartz inequality, and Lemma \ref{chap:VMS::prop:as_cont},  
\[
|(p_h, \nabla\cdot\vk_h)| =| a(\wk_h, \vk_h;\mu) + a_s'(\wk_h;\wk_h,\vk_h;\mu) + c(\wk_h,\wk_h,\vk_h;\mu) - \left<\fk, \vk_h\right>|
\]
\[
\le \dfrac{1}{\mu}\normld{\nabla\wk_h}\normld{\nabla\vk_h} + C\,\normld{\nabla\wk_h}^2\normld{\nabla\vk_h} + \normdual{\fk}\normld{\nabla\vk_h},
\]
for some constant $C>0$. Therefore, from Proposition \ref{chap:VMS::prop:infsup_b},
\[
\black{\dfrac{1}{\alpha}}\,\|p_h\|_{0,2,\Omega}\le \sup_{\vk_h\in Y_h}\dfrac{(p_h,\nabla\cdot\vk_h)_\Omega}{\normld{\nabla{\vk_h}}}+\normtau{\sigh(\nabla p_h)}
\]
\[
\le \dfrac{1}{\mu}\normld{\nabla\wk_h}+ C\,\normld{\nabla\wk_h}^2 + (1+\sqrt{\mu}\, C_P)\,\normdual{\fk}
\]
Combining this estimate with \eqref{chap:VMS::eq:estimate1}, \eqref{chap:VMS::eq:estimate3} follows.
\qed

Let us now consider the tangent operator $\partial_1A(U,V;\mu)(Z)$, defined as the Gateaux derivative of the operator $A(\cdot,\cdot;\mu)$, with respect of the first variable, in the direction $Z\in X$. From (\ref{chap:VMS::aff_decomp}) it holds
\begin{eqnarray*}
\partial_1A(U,V;\mu)(Z)&=&\frac{1}{\mu}\,\partial_1A_{Diff}(U,V)(Z)+\partial_1A_{Div}(U,V)(Z)+\partial_1A_{Conv}(U,V)(Z)\\
&+&\partial_1A_{Sma}(U;V)(Z)+ \partial_1A_{Pres}(U,V;\mu)(Z),
\end{eqnarray*}
where
\[
\begin{array}{rl}
\partial_1A_{Diff}(U,V)(Z)&\hspace{-0.25cm}=A_{Diff}(Z,V),\vspace{0.2cm}\\
\partial_1A_{Div}(U,V)(Z)&\hspace{-0.25cm}=A_{Div}(Z,V),\vspace{0.2cm}\\
\partial_1A_{Conv}(U;V)(Z)&\hspace{-0.25cm}=\displaystyle \frac{1}{2}\,\intO{\left ((\uk\cdot\nabla\zk)\vk+(\zk\cdot\nabla\uk)\vk- (\zk\cdot\nabla\vk)\uk-(\uk\cdot\nabla\vk)\zk\,\right )}, \vspace{0.2cm}\\
\partial_1A_{Sma}(U;V)(Z)&\hspace{-0.25cm}=\intO{\nu_T(\pih\wk)\,\nabla(\pih\zk):\nabla(\pih\vk)}\vspace{0.1cm} \\
	&\hspace{-1cm}+\intK{(C_Sh_K)^2\dfrac{\nabla(\pih\wk):\nabla(\pih\zk)}{|\nabla(\pih\wk)|}\big(\nabla(\pih\wk):\nabla(\pih\vk)\big)},\vspace{0.2cm}\\
\partial_1A_{Pres}(U,V;\mu)(Z)&\hspace{-0.25cm}=A_{Pres}(Z,V;\mu).
\end{array}
\]
Note that if $\nabla(\pih\wk) = 0$, $\partial_1A_{Sma}(U;V)(Z)=0$, since $A_{Sma}(U;V)=0$. We further consider that $\nabla(\pih\wk) \neq 0$. The same analysis can be analogously performed when $\nabla(\pih\wk) = 0$.

In the sequel we endow the product space $X=Y\times M$ with the hilbertian norm
\begin{equation}
\normX{U}=\sqrt{\normld{\nabla\uk}^2+\normld{p}^2}\quad \forall \,U=(\uk,p)\in X.
\end{equation}

The derivation of the a posteriori error estimation is based upon the uniform coerciveness and boundedness of the tangent operator $\partial_1A(U,V;\mu)(Z)$ on the parametric solution of problem \eqref{chap:VMS::pb:FV}: 
\begin{equation}\label{chap:VMS::ec:cont}
\infty>\gamma_0 \ge \gamma_h(\mu)\equiv\sup_{Z_h\in X_h}\sup_{V_h\in X_h}\dfrac{\partial_1A(U_h(\mu),V_h;\mu)(Z_h)}{\|Z_h\|_X\|V_h\|_X}.
\end{equation}
\begin{equation}\label{chap:VMS::ec:infsup}
0<\beta_0\le\beta_{h}(\mu)\equiv\inf_{Z_h\in X_h}\sup_{V_h\in X_h}\dfrac{\partial_1A(U_h(\mu),V_h;\mu)(Z_h)}{\|Z_h\|_X\|V_h\|_X}.
\end{equation}
for some constants $\gamma_0$, $\beta_0$. The coercivity \eqref{chap:VMS::ec:infsup} holds on branches of non-singular solutions of problem \eqref{chap:VMS::pb:FV}, while the boundedness \eqref {chap:VMS::ec:cont} is a general property:
\begin{prop}\label{chap:VMS::prop::cont}
There exists $\gamma_0\in\R$ such that $\forall \mu\in\cD$
\begin{equation}\label{coercfeda}
|\partial_1A(U_h(\mu),V_h;\mu)(Z_h)|\le \gamma_0\|Z_h\|_X\|V_h\|_X\quad \forall Z_h,V_h\in X_h.
\end{equation}
\end{prop}
\textbf{Proof.} Let us denote $U_h=(\uk_h, p_h^u)$, $V_h=(\vk_h, p_h^v)$, $Z_h=(\zk_h, p_h^z)$. It holds,
\[
|\partial_1A(U_h(\mu),V_h;\mu)(Z_h)|\le\frac{1}{\mu}|\partial_1A_{Diff}(U_h(\mu),V_h)(Z_h)|+|\partial_1A_{Div}(U_h(\mu),V_h)(Z_h)|
\]
\[
+|\partial_1A_{Conv}(U_h(\mu),V_h)(Z_h)|+|\partial_1A_{Sma}(U_h(\mu);V_h)(Z_h)|+|\partial_1A_{Pres}(U_h(\mu),V_h;\mu)(Z_h)|.
\]

The boundedness of $|\partial_1A_{Diff}(U_h(\mu),V_h;\mu)(Z_h)|$, $|\partial_1A_{Div}(U_h(\mu),V_h)(Z_h)|$, and $|\partial_1A_{Conv}(U_h(\mu);V_h)(Z_h)|$ can be proved in a standard way, the last one using estimate \eqref{esformcon} (see \cite{PaperSmago}). To bound $|\partial_1A_{Sma}(U_h(\mu);V_h)(Z_h)|$, we use the local inverse finite element estimates and the boundedness of operator $\Pi_h$. It holds,
\[
\begin{array}{ll}
|\partial_1A_{Sma}(U_h(\mu);V_h)(Z_h)|\le&\intO{(C_Sh)^2|\nabla(\pih\wk_h)||\nabla(\pih\zk_h)||\nabla(\pih\vk_h)|},\vspace{0.1cm} \\
	&\hspace{-0.7cm}+\intO{(C_Sh)^2\dfrac{|\nabla(\pih\wk_h)||\nabla(\pih\zk_h)|}{|\nabla(\pih\wk_h)|}|\nabla(\pih\wk_h)||\nabla(\pih\vk_h)|},\vspace{0.2cm}\\
	&\hspace{-0.7cm}\le 2(C_Sh)^2\normlio{\nabla(\pih\wk_h)}\normld{\nabla(\pih\zk_h)}\normld{\nabla(\pih\vk_h)}\vspace{0.1cm}\\
	&\hspace{-0.7cm}\le 2C_S^2 h^{2-d/2} \normld{\nabla(\pih\wk_h)}\normld{\nabla(\pih\zk_h)}\normld{\nabla(\pih\vk_h)}\vspace{0.1cm}\\
	&\hspace{-0.7cm}\le 2 C_1 C_S^2 h^{2-d/2}\normld{\nabla(\pih\wk_h)}\normld{\nabla\zk_h}\normld{\nabla\vk_h}\vspace{0.1cm}\\
	&\hspace{-0.7cm} \le C_3 \normX{Z_h}\normX{V_h},
\end{array}
\]
for some constants $C_1>0$, $C_2>0$, $C_3>0$, where in the last inequality we have used estimate \eqref{chap:VMS::eq:estimate1}. The last term can be bounded using Lemma \ref{chap:VMS::lema:eq_norm_tau}, and Cauchy-Schwartz inequality:
\[
\begin{array}{ll}
|\partial_1A_{Pres}(U_h(\mu);V_h;\mu)(Z_h)|\le&\intK{\tau_{p,K}(\mu)|\sigh(\nabla p_h^z)||\sigh(\nabla p_h^v)|}\vspace{0.1cm}\\
&\hspace{-0.7cm}\le\normtau{\sigh(\nabla p_h^z)}\normtau{\sigh(\nabla p_h^v)}\le C_{\tau}^2\normld{p_h^z}\normld{p^v_h}\vspace{0.1cm}\\
&\hspace{-0.7cm}\le C_4\normX{Z_h}\normX{V_h}.
\end{array}
\]
Thus, we have proved \eqref{coercfeda}.
\qed

The coercivity \eqref{chap:VMS::ec:infsup} may be proved for small enough data $\fk$:
\begin{prop}\label{chap:VMS::prop::infsup}
Suppose that \[\|\fk\|_{Y'}<\inf_{\mu \in \cal{D}}{\dfrac{1}{2\,\mu^2 \,C_P\,C_{4S}^2 }},\]
where $C_{4S}$ is the norm of the injection of $H^1_0(\Omega)^d$ in $L^4(\Omega)^d$. Then there exists $\tilde{\beta}>0$ such that,
\begin{equation} \label{coercest}
\tilde{\beta} \, \|Z_h\|_X\le \sup_{V_h\in X_h}\dfrac{\partial_1A(U_h(\mu),V_h;\mu)(Z_h)}{\|V_h\|_X}\quad \forall \, Z_h \in X_h,\quad \forall \, \mu \in {\cal D}.
\end{equation}
\end{prop}
\textbf{Proof.} Let us denote, $U_h(\mu)=(\uk_h, p_h^u)$, $Z_h=(\zk_h, p_h^z)$. It holds 
$$\partial_1A_{Pres}(U_h(\mu),Z_h;\mu)(Z_h)=\normtau{\sigh(\nabla p_h^z)}^2.
$$
 Moreover,
\[
\partial_1A_{Sma}(U_h(\mu);Z_h)(Z_h)=\intO{\nu_T(\nabla(\pih \wk_h))|\nabla (\pih \zk_h)|^2}\]\[+
\intK{(C_Sh_K)^2\dfrac{|\nabla(\pih\wk_h):\nabla(\pih\zk_h)|^2}{|\nabla(\pih\wk_h)|}}\ge0.
\]
Also,
\[
\partial_1A_{Diff}(U_h(\mu),Z_h;\mu)(Z_h)+\partial_1A_{Conv}(U_h(\mu),Z_h)(Z_h)
\]
\[\ge \dfrac{1}{\mu}\normld{\nabla \zk_h}^2+\partial_1A_{Conv}(U_h(\mu),Z_h)(Z_h)\]
\[
\ge \displaystyle \dfrac{1}{\mu}\normld{\nabla \zk_h}^2+\frac{1}{2}\,\intO{\left ((\zk_h\cdot\nabla\uk_h)\zk_h- (\zk_h\cdot\nabla\zk_h)\uk_h\,\right )}
\]
\[
\ge\frac{1}{\mu }\normld{\nabla \zk_h}^2
-\normlc{\zk_h}^2\normld{\nabla\wk_h}
-\normlc{\wk_h}\normld{\nabla\zk_h}\normlc{\zk_h}
\]
\[
\ge\left (\frac{1}{\mu } -2 \, C_{S4}^2\,\,\normld{\nabla\wk_h}\right )\, \normld{\nabla\zk_h}^2\ge\left (\frac{1}{\mu } -2 \, C_{S4}^2\,\mu\,C_P\, \|\fk\|_{Y'}\right )\, \normld{\nabla\zk_h}^2\]
\[
\ge\rho \, \normld{\nabla\zk_h}^2,\quad \forall \,\mu \in \cD,\hfill
\]
for $\rho=\displaystyle\inf_{\mu \in \cD}\left (\frac{1}{\mu } -2 \, C_{S4}^2\,\mu\,C_P\, \|\fk\|_{Y'}\right ) >0$. Consequently,
$$
\partial_1A(U_h(\mu),Z_h;\mu)(Z_h)\ge\rho \, \normld{\nabla\zk_h}^2+\normtau{\sigh(\nabla p_h^z)}^2.
$$
Let 
$$
S(\mu)=\displaystyle \sup_{V_h\in X_h}\dfrac{\partial_1A(U_h(\mu),V_h;\mu)(Z_h)}{\|V_h\|_X}.
$$
 Then, in abridged notation for brevity,
\begin{equation} \label{estsup1}
\rho \, \|\zk_h\|_Y^2+\normtau{\sigh(\nabla p_h^z)}^2 \le S(\mu)\, \|Z_h\|_X= S(\mu) \, \left (\|\zk_h\|_Y^2+\|p_h^z\|_M^2\right )^{1/2}.
\end{equation}
Also, letting $V_h=(\vk_h,0)\in X_h$, estimates similar to those carried on in the proof of Proposition \ref{chap:VMS::prop::cont} yield
$$
(\nabla\cdot \vk_h,p_h^z) \le \left ( S(\mu)+(\frac{1}{\mu}+C)\,\|\zk_h\|_Y\right )\, \|\vk_h\|_Y
$$
for some constant $C>0$ independent of $\nu$ and $h$. Then, letting \break $\displaystyle z= \sup_{\vk_h\in Y_h}\dfrac{(p_h^z,\nabla\cdot \vk_h)_\Omega}{\|\vk_h\|_Y}$,
\begin{equation} \label{estsup2}
z \le  S(\mu)+C'(\mu)\,\|\zk_h\|_Y\,\quad \mbox{with   } C'(\mu)=\frac{1}{\mu}+C.
\end{equation}
Combining \eqref{estsup1} and \eqref{estsup2},  and using \eqref{eq:infsup_b},
\begin{eqnarray*} 
&&\rho \, \|\zk_h\|_Y^2+\normtau{\sigh(\nabla p_h^z)}^2 \le  S(\mu) \, \left (\|\zk_h\|_Y+z+ \normtau{\sigh(\nabla p_h^z)} \right ) \nonumber \\
&&\le S(\mu) \, \left (\alpha\,S(\mu)+(1+\alpha\,C'(\mu))\, \|\zk_h\|_Y^2+ \alpha\,\normtau{\sigh(\nabla p_h^z)} \right ) \nonumber 
\end{eqnarray*}
Using Young's inequality, it follows
\begin{eqnarray} 
\rho \, \|\zk_h\|_Y^2+\normtau{\sigh(\nabla p_h^z)}^2 &\le&  C\, \left ( 1+\frac{1+C'(\mu)}{\rho}\right )\, S(\mu)^2\nonumber \\
&\le&  C\, \left ( 1+\frac{1}{\mu\, \rho}\right )\, S(\mu)^2. \label{estsup3}
\end{eqnarray}
Combining this inequality with \eqref{estsup2},
\begin{equation} \label{estsup4}
z \le  \left (1+C\, C'(\mu)\,\sqrt{1+\frac{1}{\mu\, \rho}}\right )\, S(\mu).
\end{equation}
Taking the square rooth of \eqref{estsup3} and summing it with \eqref{estsup4} finally yields
\begin{eqnarray} 
&\displaystyle  \min\{\sqrt{\rho},\,\alpha^{-1}\} \,\|Z_h\|_X \le \sqrt{\rho} \, \|\zk_h\|_Y+z+ \normtau{\sigh(\nabla p_h^z)} &\nonumber \\
&\le  \displaystyle C\, \left (1+\frac{1}{\mu} \right )\,\sqrt{1+\frac{1}{\mu\, \rho}}\, S(\mu)
\le M\, S(\mu), \,\, \label{estsup5}\nonumber&
\end{eqnarray}
with   $M= \displaystyle \sup \left \{ C\, \left (1+\frac{1}{\mu} \right )\,\sqrt{1+\frac{1}{\mu\, \rho}},\,\, \,\mu \in \cD\right \}$. Then, \eqref{coercest} holds with $\tilde{\beta}=\min\{\sqrt{\rho},\,\alpha^{-1}\}/M$.
\qed

The same demonstration proves the coercivity of the reduced problem for small data, when the reduced space $X_N$ is endowed with the norm
\begin{equation}
\normXN{U_N}=\sqrt{\normld{\nabla\uk_N}^2+\normp{p_N}^2}\quad \forall \,U_N=(\uk_N,p_N)\in X_N,
\end{equation}
where $\normp{\cdot}$ is defined in \eqref{infsupred}.

\begin{prop}\label{chap:VMS::prop::infsupN}
Assume that the pressure finite element space $M_h$ is formed by piecewise affine functions. Suppose that \[\|\fk\|_{Y'}\le \inf_{\mu \in \cal{D}}{\dfrac{1}{2\,\mu^2 \,C_P\,C_{4S}^2 }}.\]
Then there exists a constant $\tilde{\beta}>0$ such that,
\begin{equation} \label{coercest}
\tilde{\beta} \, \normXN{Z_N}\le \sup_{V_N\in X_N}\dfrac{\partial_1A(U_N(\mu),V_N;\mu)(Z_N)}{\normXN{V_N}}\quad \forall \, Z_N \in X_N,\quad \forall \, \mu \in {\cal D}.
\end{equation}
\end{prop}
In this case, the $\normXN{\cdot}$ norm is bounded from above by a constant times the $\normld{\cdot}$ norm but, to the best of our knowledge, they are no longer equivalent.
 \section{\textit{A posteriori} error bound estimator}\label{chap:VMS::sec:Post_Err}
In this section we develop the \textit{a posteriori} error bound estimator used in the hierarchical construction of the reduced space during the Greedy algorithm. We use the Brezzi-Rappaz-Raviart (BRR) theory \cite{BRR} for approximation on a regular branch of solutions of non-linear problems. We start by proving the local Lipschitz continuity of the tangent Smagorinsky operator, which is a key result to apply this theory:

\begin{lema}\label{chap:VMS::lema:LemmaRho}
It holds 
\begin{eqnarray}\label{chap:VMS::ec:ro}
\left|\partial_1A(U_h^1,V_h;\mu)(Z_h)-\partial_1A(U_h^2,V_h;\mu)(Z_h)\right|\le&& \\
\rho_T\normX{U_h^1-U_h^2}\normX{Z_h}\|V_h\|_X,
\nonumber \end{eqnarray}
for all $ U_h^1,U_h^2,Z_h,V_h \in X_h$, with $\rho_T=2\,C_{S4}^2+C\,h^{2-d/2}$, for some constant $C>0$.
\end{lema}
\textbf{Proof.} We have
\begin{align*}
&\partial_1A(U_h^1,V_h;\mu)(Z_h)-\partial_1A(U_h^2,V_h;\mu)(Z_h)= \\
&\displaystyle\frac{1}{2}\left (\intO{((\uk_h^1-\uk_h^2)\cdot\nabla \zk_h)\vk_h}+\intO{(\zk_h\cdot\nabla(\uk_h^1-\uk_h^2))\vk_h}\right .\\
&\left .-\intO{(\zk_h\cdot\nabla \vk_h)(\uk_h^1-\uk_h^2)}-\intO{((\uk_h^1-\uk_h^2)\cdot\nabla\vk_h)\zk_h}\right )\\
&+\intO{\left[\nu_T(\pih\uk_h^1)-\nu_T(\pih\uk_h^2)\right]\nabla(\pih\zk_h):\nabla(\pih\vk_h)}\\
&+\intK{(C_Sh_K)^2\dfrac{\nabla(\pih\uk_h^1):\nabla(\pih\zk_h)}{|\nabla(\pih\uk_h^1)|}\big(\nabla(\pih\uk_h^1):\nabla(\pih\vk_h)\big)}\\
&-\intK{(C_Sh_K)^2\dfrac{\nabla(\pih\uk_h^2):\nabla(\pih\zk_h)}{|\nabla(\pih\uk_h^2)|}\big(\nabla(\pih\uk_h^2):\nabla(\pih\vk_h)\big)}.
\end{align*}

Thanks to the triangular inequality, it holds
\begin{equation}\label{chap:VMS::eq:derSmago}
\begin{array}{lc}
&\left|\partial_1A(U_h^1,V_h;\mu)(Z_h)-\partial_1A(U_h^2,V_h;\mu)(Z_h)\right|\le\\ 
&\displaystyle\frac{1}{2}\left|\intO{((\uk_h^1-\uk_h^2)\cdot\nabla \zk_h)\vk_h}\right|+\displaystyle\frac{1}{2}\left|\intO{(\zk_h\cdot\nabla(\uk_h^1-\uk_h^2))\vk_h}\right|\\
&+\displaystyle\frac{1}{2}\left|\intO{(\zk_h\cdot\nabla \vk_h)(\uk_h^1-\uk_h^2)}\right |+ \displaystyle\frac{1}{2}\left|\intO{((\uk_h^1-\uk_h^2)\cdot\nabla\vk_h)\zk_h}\right |\\
&+\left| \intO{\left[\nu_T(\pih\uk_h^1)-\nu_T(\pih\uk_h^2)\right]\nabla(\pih\zk_h):\nabla(\pih\vk_h)}\right|\\
&+\left|\intK{(C_Sh_K)^2\dfrac{\nabla(\pih\uk_h^1):\nabla(\pih\zk_h)}{|\nabla(\pih\uk_h^1)|}\big(\nabla(\pih\uk_h^1):\nabla(\pih\vk_h)\big)}\right.\\
&\left.-\intK{(C_Sh_K)^2\dfrac{\nabla(\pih\uk_h^2):\nabla(\pih\zk_h)}{|\nabla(\pih\uk_h^2)|}\big(\nabla(\pih\uk_h^2):\nabla(\pih\vk_h)\big)}\right|.\\
\end{array}
\end{equation}

We bound each term in (\ref{chap:VMS::eq:derSmago}) separately. The four first terms can be bounded in a standard way using estimate \eqref{esformcon}. For the third term, we use the local inverse inequalities of finite element functions (cf. \cite{desinv}) and the stability of operator $\Pi_h$,
\[
\left| \intO{\left(\nu_T(\pih\uk_h^1)-\nu_T(\pih\uk_h^2)\right)\nabla(\pih\zk_h):\nabla(\pih\vk_h)}\right|
\]
\[
\le \intK{(C_Sh_K)^2\big||\nabla (\pih \uk_h^1)| - |\nabla (\pih \uk_h^2)| \big||\nabla (\pih \zk_h)||\nabla (\pih \vk_h)|}
\]
\[
\le (C_Sh)^2\intO{|\nabla \big(\pih (\uk_h^1-\uk_h^2)\big)||\nabla (\pih \zk_h)||\nabla (\pih \vk_h)|}
\]
\[
\le (C_Sh)^2\normlio{\nabla \big(\pih (\uk_h^1-\uk_h^2)\big)}\normld{\nabla (\pih \zk_h)}\normld{\nabla (\pih \vk_h)}
\]
\[
\le C_S^2h^{2-d/2}C\normld{\nabla \big(\pih (\uk_h^1-\uk_h^2)\big)}\normld{\nabla (\pih \zk_h)}\normld{\nabla (\pih \vk_h)}
\]
\[
\le C_S^2h^{2-d/2}CC_f \,\normld{\nabla (\uk_h^1-\uk_h^2)}\normld{\nabla \zk_h}\normld{\nabla  \vk_h}
\]
\[
\le C_S^2h^{2-d/2}CC_f \,\normX{U_h^1-U_h^2}\normX{Z_h}\normX{V_h}
\]

The last term in (\ref{chap:VMS::eq:derSmago}) is bounded as follows:
\[
\left|\intK{(C_Sh_K)^2\dfrac{\nabla(\pih\uk_h^1):\nabla(\pih\zk_h)}{|\nabla(\pih\uk_h^1)|}\big(\nabla(\pih\uk_h^1):\nabla(\pih\vk_h)\big)}\right.
\]
\[
\left.-\intK{(C_Sh_K)^2\dfrac{\nabla(\pih\uk_h^2):\nabla(\pih\zk_h)}{|\nabla(\pih\uk_h^2)|}\big(\nabla(\pih\uk_h^2):\nabla(\pih\vk_h)\big)}\right|
\]
\[
=\left|\intK{(C_Sh_K)^2\left[\frac{\nabla(\pih\uk_h^1):\nabla(\pih\zk_h)}{|\nabla(\pih\uk_h^1)|}\left(\nabla(\difupih):\nabla(\pih\vk_h)\right)\right.\right.
\]
\[
+\frac{\nabla(\difupih):\nabla(\pih\zk_h)}{|\nabla(\pih\uk_h^2)|}\big(\nabla(\pih\uk_h^2):\nabla(\pih\vk_h)\big)
\]
\[
+\left.\left.\frac{(|\nabla(\pih\uk_h^2)|-|\nabla(\pih\uk^1_h)|)\nabla(\pih\uk_h^1):\nabla(\pih\zk_h)}{|\nabla(\pih\uk_h^1)||\nabla(\pih\uk^2_h)|}\big(\nabla(\pih\uk_h^2):\nabla(\pih\vk_h)\big)\right]}\right|
\]
\[
\le\intK{(C_Sh_K)^2|\nabla(\pih\zk_h)||\nabla(\difupih)||\nabla(\pih\vk_h)|}
\]
\[
+\intK{(C_Sh_K)^2|\nabla(\difupih)||\nabla(\pih\zk_h)||\nabla(\pih\vk_h)|}
\]
\[
+\intK{(C_Sh_K)^2\left||\nabla(\pih\uk_h^1)|-|\nabla(\pih\uk^2_h)|\right||\nabla(\pih\zk_h)||\nabla(\pih\vk_h)|}
\]
\[
\le3(C_Sh)^2\normlio{\nabla(\difupih)}\normld{\nabla(\pih\zk_h)}\normld{\nabla(\pih\vk_h)}
\]
\[
\le3C_S^2h^{2-d/2}C\normld{\nabla(\difupih)}\normld{\nabla(\pih\zk_h)}\normld{\nabla(\pih\vk_h)}\]
\[
\le3C_S^2h^{2-d/2}CC_f \,\normld{\nabla(\difuh)}\normld{\nabla\zk_h}\normld{\nabla\vk_h}\]
\[
\le3C_S^2h^{2-d/2}CC_f \,\normX{\difUh}\normX{Z_h}\normX{V_h}.
\]

Thus, we have just proved that
\[
\left|\partial_1A(U_h^1,V_h;\mu)(Z_h)-\partial_1A(U_h^2,V_h;\mu)(Z_h)\right|\le
\rho_T\normX{\difUh}\normX{Z_h}\normX{V_h},
\]
with, $\rho_T=2C_T+C\, h^{2-d/2}$ for some constant $C>0$.
\qed

We assume that the RB problem (\ref{chap:VMS::pb:RBP}) is well posed, and define the fo\-llo\-wing inf-sup and continuity constants:
\begin{equation}\label{chap:LPS::betaN}
0<\beta_{N}(\mu)\equiv\inf_{Z_h\in X_h}\sup_{V_h\in X_h}\dfrac{\partial_1A(\unmu,V_h;\mu)(Z_h)}{\|Z_h\|_X\|V_h\|_X},
\end{equation}
\begin{equation}\label{chap:LPS::gammaN}
\infty>\gamma_{N}(\mu)\equiv\sup_{Z_h\in X_h}\sup_{V_h\in X_h}\dfrac{\partial_1A(\unmu,V_h;\mu)(Z_h)}{\|Z_h\|_X\|V_h\|_X},
\end{equation}
Let us introduce the \textit{a posteriori} error bound estimator:
\begin{equation}\label{chap:LPS::delta}
\dn=\frac{\beta_N(\mu)}{2\rho_T}\left[1-\sqrt{1-\taun}\right],
\end{equation}
where $\taun$ is given by
\begin{align}
\taun&=\frac{4\en\rho_T}{\beta_N^2(\mu)},\label{chap:LPS::tau}
\end{align} 
with $\en$ the dual norm, on $X_h'$, of the residual ${\cal R}(U_N(\mu);\mu)$, defined as
$$
{\cal R}( U_N(\mu);\mu)(V_h)=A(U_N(\mu),V_h;\mu)-F(V_h;\mu);
$$
that is,
$$
\en =\sup_{V_h \in X_h}\frac{{\cal R}( U_N(\mu);\mu)(V_h)}{\|V_h\|_X}.
$$

The \textit{a posteriori} error bound estimator is stated by the following result

\begin{teor}\label{chap:VMS::teor:Teorprinc}
Let $\mu\in\cD$. \black{Assume that the tangent Smagorinsky operator satisfies \eqref{chap:VMS::ec:ro}, and also assume that is continuous and inf-sup stable, verifying that $\beta_N(\mu)>0$}. If problem (\ref{chap:VMS::pb:FVX}) admits a solution $U_h(\mu)$ such that
\[
\normX{U_h(\mu)-U_N(\mu)}\le\frac{\beta_N(\mu)}{\rho_T},
\]
then this solution is unique in the ball $B_X\left(U_N(\mu),\dfrac{\beta_N(\mu)}{\rho_T}\right)$. 

Moreover, assume that $\taun\le1$ for all $\mu\in\cD$. Then there exists a unique solution $U_h(\mu)$ of (\ref{chap:VMS::pb:FVX}) such that the error with respect to $U_N(\mu)$, solution of (\ref{chap:VMS::pb:RBP}), is bounded by 
\begin{equation}\label{chap:VMS::teorprinc:err}
\normX{U_h(\mu)-\unmu}\le\dn\le\left[\frac{2\gamma_{N}(\mu)}{\beta_N(\mu)}+\taun\right]\normX{U_h(\mu)-\unmu},
\end{equation}
\black{where $\gamma_N(\mu)$ denotes the dual norm of $\partial_1A(U_N(\mu);\cdot;\mu)$ given by \eqref{chap:LPS::gammaN}.}
\end{teor}

The proof of this theorem follows that of Theorem 5.3 in \cite{PaperSmago}, we omit it for brevity.

\begin{nota}
Note that $\dn$ is defined only if $\taun\le1$. In practice, at the first iterations of the Greedy algorithm, $\taun>1$. For these cases, we consider a Taylor approximation of the \textit{a posteriori} error estimator, that is given by
\begin{equation}
\dn\simeq\dfrac{\beta_N(\mu)}{4\rho_T}\taun.
\end{equation}
In Section \ref{chap:VMS::sec:num_results} we show that this Taylor approximation is more accurate than $\taun$ as estimator.
\end{nota}


\section{Approximation of eddy viscosity term and pressure stabilizing coefficient}\label{chap:VMS::sec:EIM}

In this section, we present the approximation of the non-linear terms with respect to the parameter, throughout the Empirical Interpolation Method \cite{EIM1,EIM2}. 

Both the Small-Small setting of the Smagorinsky eddy-diffusion term defined in (\ref{chap:VMS::aff_decomp}), $\nu_T(\pih(\nabla\wk)):=\nu_T(\mu)$, and the pressure stabilisation coefficient, $\tau_{p,K}(\mu)$, defined in (\ref{chap:VMS::eq:s_pres}) are non-linear functions of the parameter, and consequently need to be linearised with the EIM to build the RB model to reduce the on-line computation times.  

For this purpose, we build two reduced-basis spaces $W_{M_1}^S=\{q_1^S(\mu),\dots,q_{M_1}^S(\mu)\}$ and  $W_{M_2}^P=\{q_1^P(\mu),\dots,q_{M_2}^P(\mu)\}$ by a greedy procedure selection, with $W_{M_1}^S$ and $W_{M_2}^P$ the EIM reduced spaces associated to the eddy-viscosity term and the pressure-stabilisation term, respectively.

Thus, we approximate them by the following trilinear forms:
\begin{equation}
\begin{array}{rl}
a_S'(\wk_N;\wk_N,\vk_N;\mu)&\hspace{-0.15cm}\approx \hat{a}_S'(\wk_N,\vk_N;\mu),\vspace{0.1cm}\\
s_{pres}(p_N,q_N;\mu)&\hspace{-0.15cm}\approx \hat{s}_{pres}(p_N,q_N;\mu),
\end{array}
\end{equation}
where,
\begin{equation}
\begin{array}{rl}
\hat{a}_S'(\wk_N,\vk_N;\mu)&\hspace{-0.15cm}=\dsum_{k=1}^{M_1}\sigma^S_k(\mu)s(q_k^S,\wk_N,\vk_N),\vspace{0.1cm}\\
\hat{s}_{pres}(p_N,q_N;\mu)&\hspace{-0.15cm}=\dsum_{k=1}^{M_2}\sigma^P_k(\mu)r(q_k^P,p_N,q_N),
\end{array}
\end{equation}
with,
\begin{equation}
\begin{array}{rl}
s(q_k^S,\wk_N,\vk_N)&=\dsum_{K\in\Th}\big(q^S_k \,\nabla(\pih\wk_N),\nabla(\pih\vk_N)\big)_K,\vspace{0.1cm}\\
r(q_k^P,p_N,q_N)&=\dsum_{K\in\Th}\big(q^P_k\,\sigh(\nabla p_N),\sigh(\nabla q_N)\big)_K. 
\end{array}
\end{equation}

Here we are considering that the approximations given by the EIM for $\nu_T(\mu)$ and $\tau_{K,p}(\mu)$ are respectively
\begin{equation} \label{redcoefs}
\begin{array}{rl}
\mathcal{I}_{M_1}[\nu_T(\mu)]&\hspace{-0.15cm}=\dsum_{k=1}^{M_1}\sigma^S_k(\mu)q_k^S,\vspace{0.1cm}\\
\mathcal{I}_{M_2}[\tau_{K,p}(\mu)]&\hspace{-0.15cm}=\dsum_{k=1}^{M_2}\sigma^P_k(\mu)q_k^P.
\end{array}
\end{equation} 
In practice, to deal with non-homogeneous Dirichlet boundary conditions, we split the solution as $\uk_N + \uk_D$, where $\uk_N \in X_h$ and $\uk_D \in H^1(\Omega)^d$ is a lift of the Dirichlet boundary condition, compatible with the remaining boundary conditions (if any). We are thus led to a \lq\lq practical" RB problem, given by
\begin{equation}\label{chap:VMS::pb:RB_LPS_LIN}\left\{\begin{array}{l}
\mbox{Find } (\uk_N,p_N)=(\uk_N(\mu),p_N(\mu))\in Y_N\times M_N\mbox{ such that}\vspace{0.3cm}\\

\begin{array}{ll}
a(\uk_N,\vk_N;\mu)+b(\vk_N,p_N;\mu)+ \hat{a}_S'(\wk_N;\vk_N;\mu) \\
+c(\uk_N,\uk_N,\vk_N;\mu)+c(\uk_D,\uk_N,\vk_N;\mu)\\
+c(\uk_N,\uk_D,\vk_N;\mu)=F(\vk_N;\mu)&\quad\forall\vk_N\in Y_N\vspace{0.1cm}\\
b(\uk_N,q_N;\mu)+\hat{s}_{pres}(p_N,q_N;\mu)=0&\quad\forall q_N\in M_N.\end{array}\end{array}\right.
\end{equation} 

We can express the solution $(\uk_N(\mu),p_N(\mu))\in X_N$ of (\ref{chap:VMS::pb:RB_LPS_LIN}) as a linear combination of the basis functions:

\[
\uk_N(\mu)=\sum_{j=1}^{N}u_j^N(\mu)\zeta^{\vk}_j,\quad p_N(\mu)=\sum_{j=1}^Np_j^N(\mu)\xi^p_j.
\]

The tensor representation associated to the eddy-viscosity term and the pressure stabilisation coefficient  are define as follow:
\begin{equation}\nonumber
(\mathbb{S}'_N(q^S_s))_{ij}=\sum_{K\in\Th}(q^S_s \nabla(\pih \zetav_j),\nabla(\pih \zetav_i))_K,\quad i,j=1\dots,N, s=1,\dots,M_1,
\end{equation}
\begin{equation}\nonumber
(\mathbb{P}_N(q^P_s))_{ij}=\sum_{K\in\Th}(q^P_s\sigh(\nabla \xip_j),\sigh(\nabla \xip_i))_K,\quad i,j=1,\dots,N, s=1,\dots,M_2.
\end{equation}
With this tensor representation for the non-linear terms in (\ref{chap:VMS::pb:RBP}) done in the offline phase, it holds that
\[
\hat{a}_S'(\zetav_j;\zetav_i;\mu)=\sum_{s=1}^{M_1}\sigma_s^S(\mu)\mathbb{S}'_N(q^S_s)\quad \mbox{ and } \quad
\hat{s}_{pres}(\xip_j,\xip_i;\mu)=\sum_{s=1}^{M_2}\sigma_s^P(\mu)\mathbb{P}_N(q^P_s).
\]
Problem (\ref{chap:VMS::pb:RB_LPS_LIN}) is solved by a semi-implicit evolution approach. \black{The terms \;
$c(\uk,\uk,\vk;\mu)$ and $a_S'(\uk,\uk,\vk;\mu)$ are linearised by
\[
c(\uk_N^n, \uk_N^{n+1}, \vk_N;\mu), \quad a_S'(\uk_N^n;\uk_N^{n+1},\vk;\mu),
\]
whereas the discretisation of the remaining terms is implicit.
}

Thanks to these linearised representations, we are able to efficiently solve the RB model in the online phase.

\section{Numerical Results}\label{chap:VMS::sec:num_results}

In this section, we present some numerical results to test the LPS pressure stabilisation for the Smagorinsky reduced basis model. We consider  the lid-driven cavity problem for two different finite element discretisations. To let this problem fit into the preceding theory we consider a lifting $\uk_D \in [H^1(\Omega)]^2$ of the non-homogeneous boundary conditions, and search for the solution as $\uk=\uk_0+\uk_D$, where the new unknown $\uk_0$ satisfies homogeneous Dirichlet boundary conditions. The change in the convection term makes it appear a linear term with the structure $\uk_D\nabla \cdot \uk_0 + \uk_0\nabla \cdot \uk_D$, the remaining operator terms in the Smagorinsky model have the same structure. An involved, but straightforward, extension of the preceding analysis, proves that it also applies to this situation whenever $\uk_D$ is small enough with respect to the viscosity coefficient. We have preferred not to include all this analysis for brevity. \black{In our numerical tests we use as $\uk_D$ a lift belonging to $[H^1(\Omega)]^2$ of a smoothing of the piecewise constant Dirichlet bounday condition given by
\[
\uk_D = \left\{\begin{array}{ll}
(1,0) & \text{if }\xk_h\in\Gamma_{top}, \smallskip\\
(0,0) & \text{if }\xk_h\notin\Gamma_{top},
\end{array} \right.
\]
where we denote by $\Gamma_{top}$ the lid of the cavity.}

We include LPS stabilisation of the pressure for both the full order and the reduced methods. On one hand the Taylor-Hood ($\mathbb{P}_2-\mathbb{P}_1$) finite element, that fits into the theory developed above, then both the full order finite element and the reduced methods are stable. On another hand, although the LPS-stabilised full-order method constructed with the $\mathbb{P}_2-\mathbb{P}_2$ finite element for velocity-pressure is stable, we do not know any proof that the LPS-stabilised reduced basis problem without pressure supremizers, constructed from this stabilized $\mathbb{P}_2-\mathbb{P}_2$ full order method, is stable. For both cases, we consider the two strategies of either including pressure supremizers, or not, to build the reduced pressure spaces, to compare the performances of both strategies. We use FreeFem++ to perform the computations (cf. \cite{freefem++}).

We consider the Reynolds number as a parameter, ranging in \break $\cD=[1000,5100]$. In this range of Reynolds number, the solution reaches a steady state regime. We consider a regular mesh with 5000 triangles and 2601 nodes. In both cases, we implement the EIM for the eddy viscosity term and the pressure stabilisation constant. 

We are aware that in the turbulent regime the lid-driven cavity flow is 3D (Cf. \cite{TomasSmago}, Chapter 13). We consider these 2D tests to assess the theoretical derivations within the paper, and to compare the two strategies of including - not including pressure supremizers in the reduced velocity spaces, when the reduced model includes the LPS stabilisation of the pressure.

For the numerical tests, we consider a weighted the scalar product for the velocity space, that lets us improve the efficiency of the \textit{a posteriori }error estimations as shown in \cite{PaperSmago}. This norm is defined as 
\begin{equation}\label{chap:VMS::normT}
(\uk_h,\vk_h)_T=	\intO{\left[\frac{1}{\mub}+\nu_T'^*\right]\nabla\uk_h:\nabla\vk_h}\qquad\forall\uk_h,\vk_h\in Y_h,
\end{equation} 
where $\nu_T'^*=\nu_T(\pih\wk_h(\mub))$, and 
\[\mub=\arg\D\min_{\mu\in\cD}\D\sum_{K\in\mathcal{T}_h}(C_Sh_K)^2\D\min_{x\in K}|\nabla(\pih\wk_h)(\mu)|(x)\chi_K(x),\] 
with $\wk_h(\mu)$ the velocity solution of (\ref{chap:VMS::pb:FV})
\subsection{Test 1. Taylor-Hood finite element}

For this test, we need $M_1=49$ basis functions $q_k^S$ for the eddy viscosity $\nu_T$ EIM approximation and $M_2=29$ basis functions $q_l^S$ for the stabilisation pressure coefficient $\tau_{p,K}$ EIM approximation to reach a prescribed tolerance of $\veps_{EIM}=5\cdot10^{-5}$ for both terms. In Fig. \ref{chap:VMS::fig:EIMCav_TH} we show the evolution of the error for the both EIM approximations. 
\begin{figure}[h]
\centering
\includegraphics[width=0.45\linewidth]{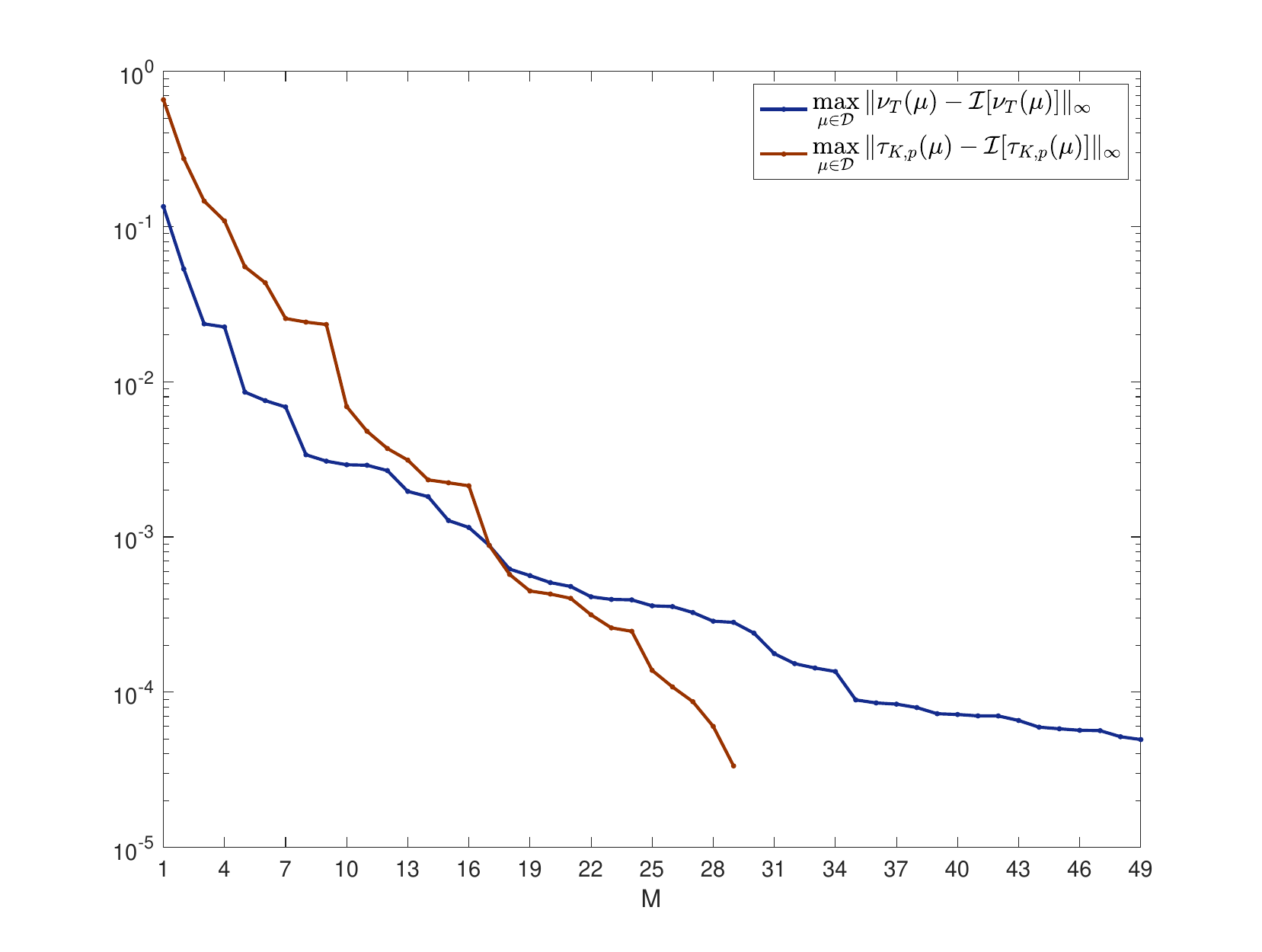}
\caption{Test 1. Evolution of the parametric $L^\infty$ error in the EIM approximation of eddy viscosity and stabilisation pressure coefficient.}\label{chap:VMS::fig:EIMCav_TH}
\end{figure}

For the Greedy algorithm that selects the reduced basis functions for velocity and pressure, we prescribe a tolerance of $\veps_{RB}=7\cdot10^{-5}$, which is reached for $N=16$ basis functions in both cases. In Fig. \ref{chap:VMS::fig:EstErrCav_TH} we show the evolution of the \textit{a posteriori} error bound estimator during the Greedy algorithm for both strategies (including pressure supremizers or not in the reduced velocity space). We can observe that the evolution of the maximum value for the error estimator either considering supremizers (left) or not (right) are equivalent, giving the same number of basis functions until reaching the tolerance. 

\begin{figure}[h]
\centering
\includegraphics[width=0.45\linewidth]{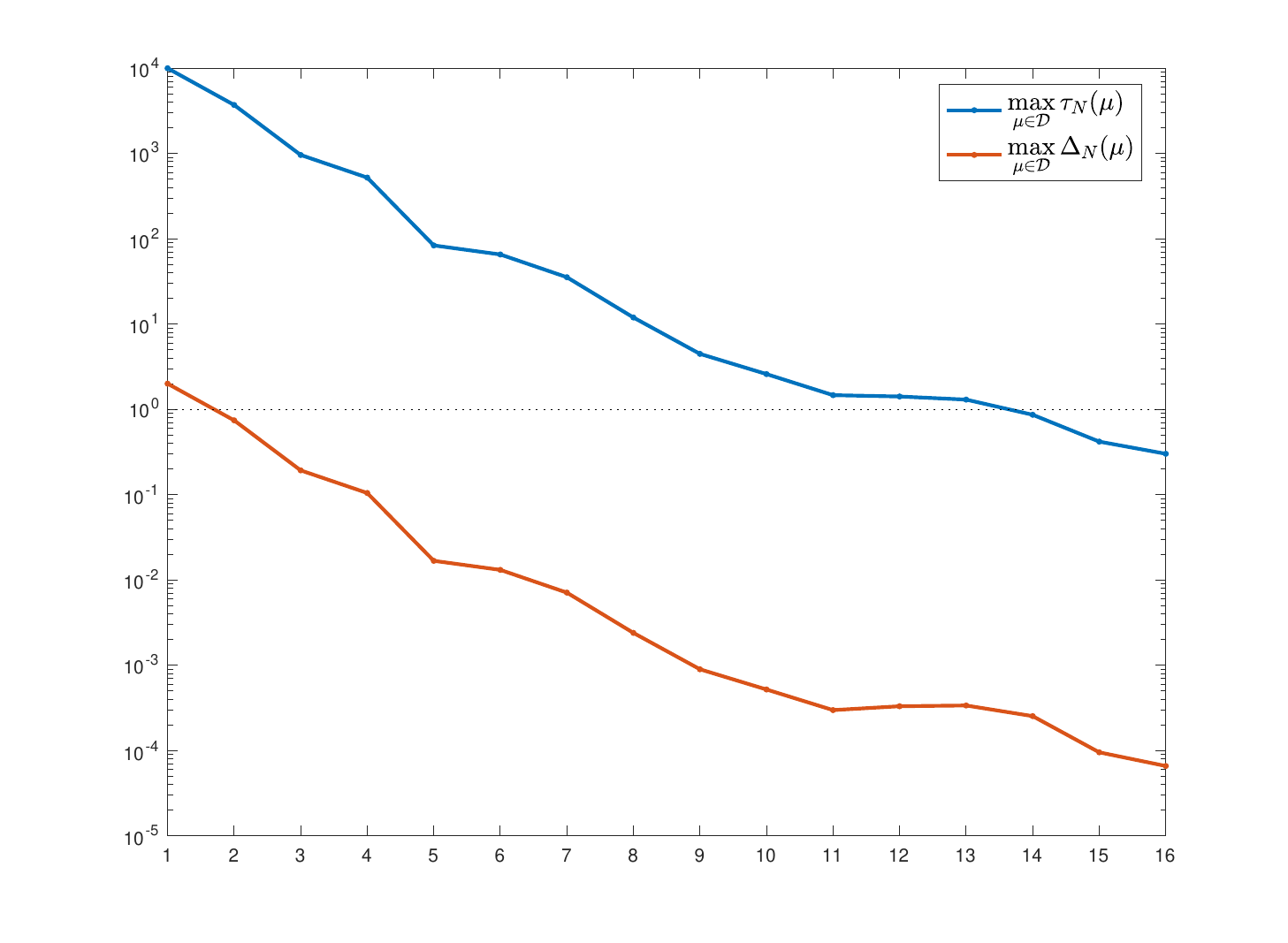}
\includegraphics[width=0.45\linewidth]{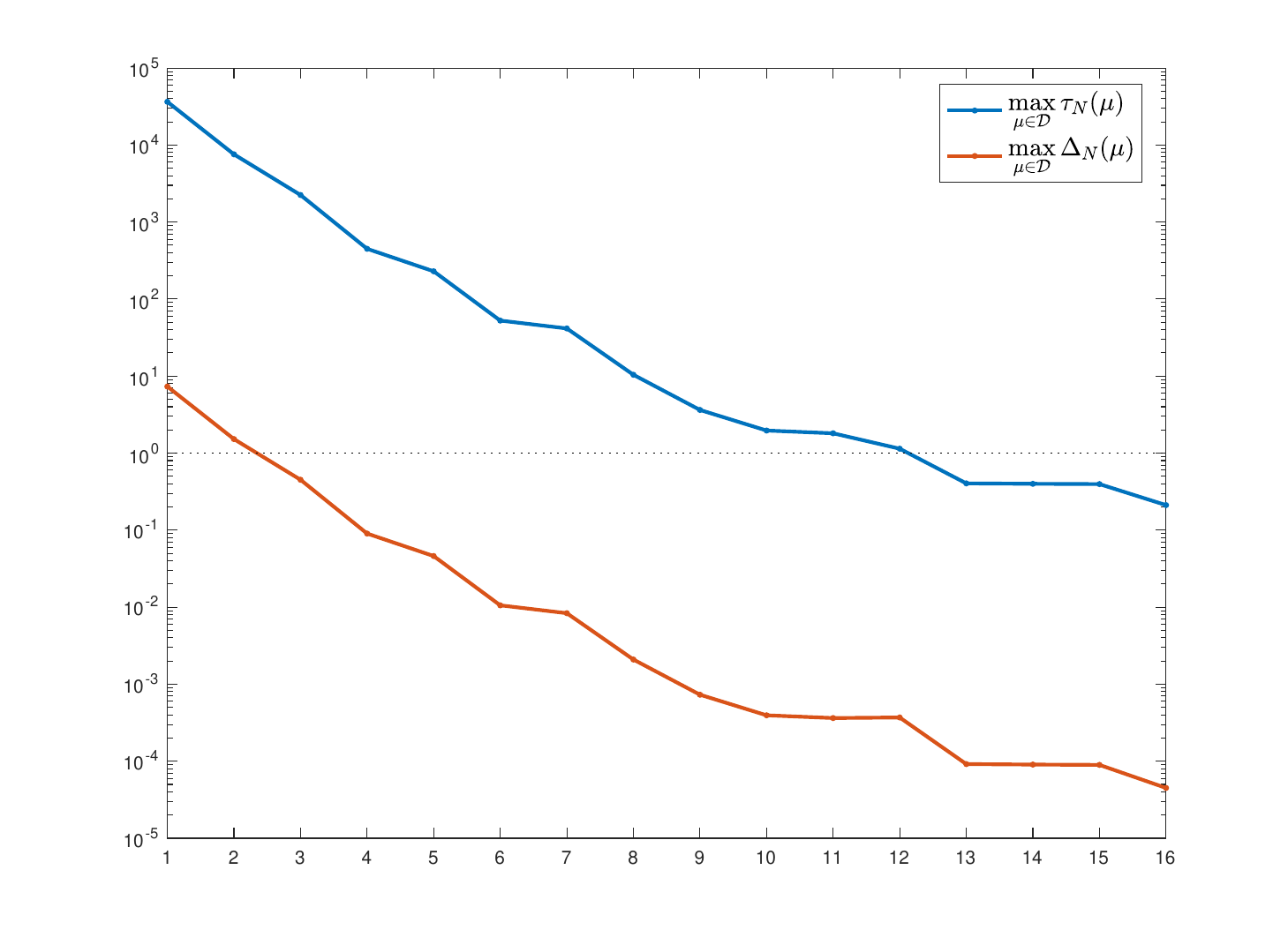}
\caption{Test 1. Evolution of the error in the Greedy algorithm with (left) and without pressure supremizers (right).}\label{chap:VMS::fig:EstErrCav_TH}
\end{figure}

We also have measured, for some randomly selected parameter values, the maximum exact error for velocity and pressure snapshot solutions with respect to the reduced solution computed, during the Greedy algorithm. We can observe in Fig. \ref{chap:LPS::fig:error_vel_pres_TH} that for both the velocity and pressure the maximum error values in each iteration of the Greedy algorithm, for both velocity and pressure, are quite close. 

\begin{figure}[h]
\centering
\includegraphics[width=0.45\linewidth]{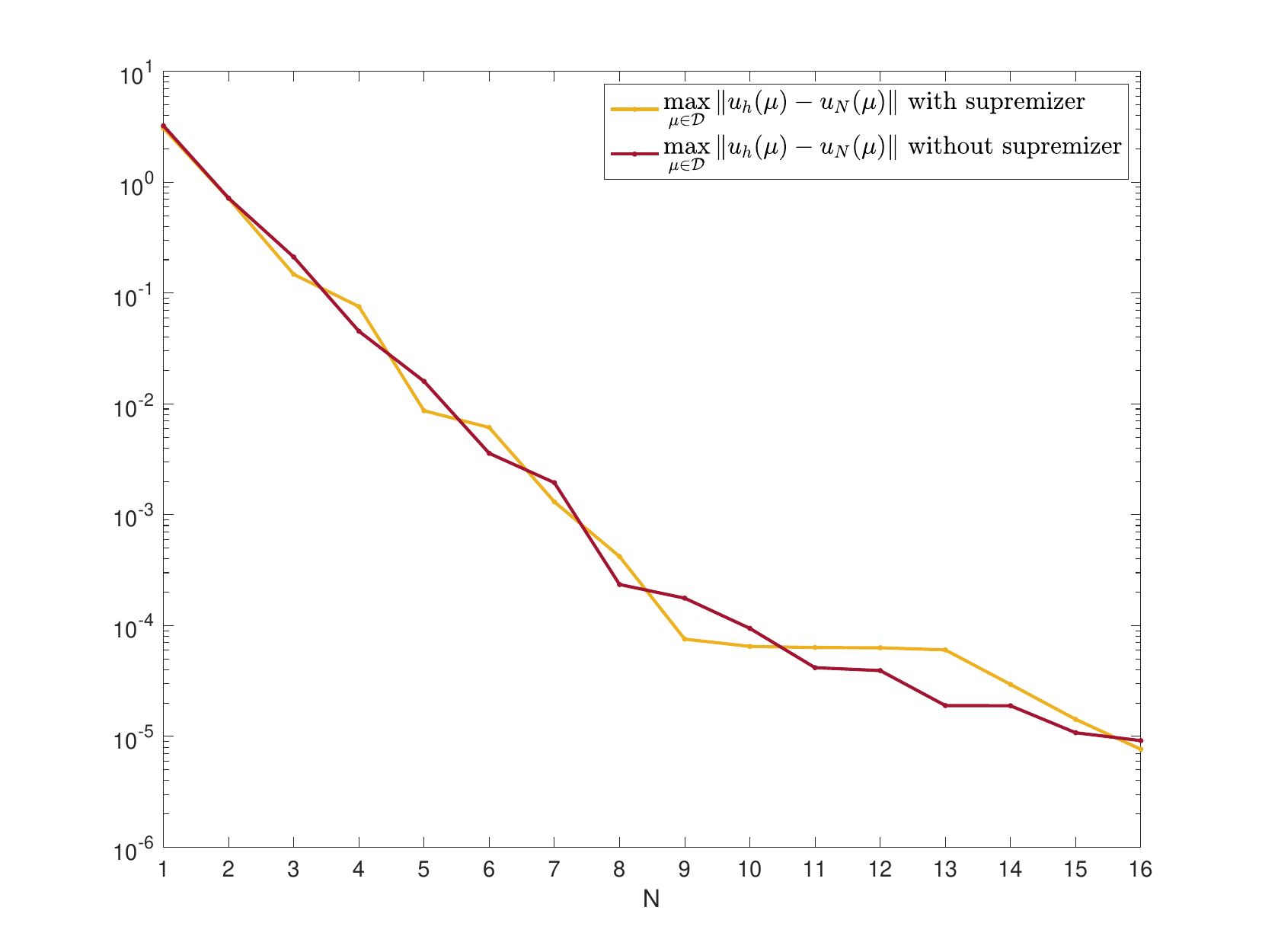}
\includegraphics[width=0.45\linewidth]{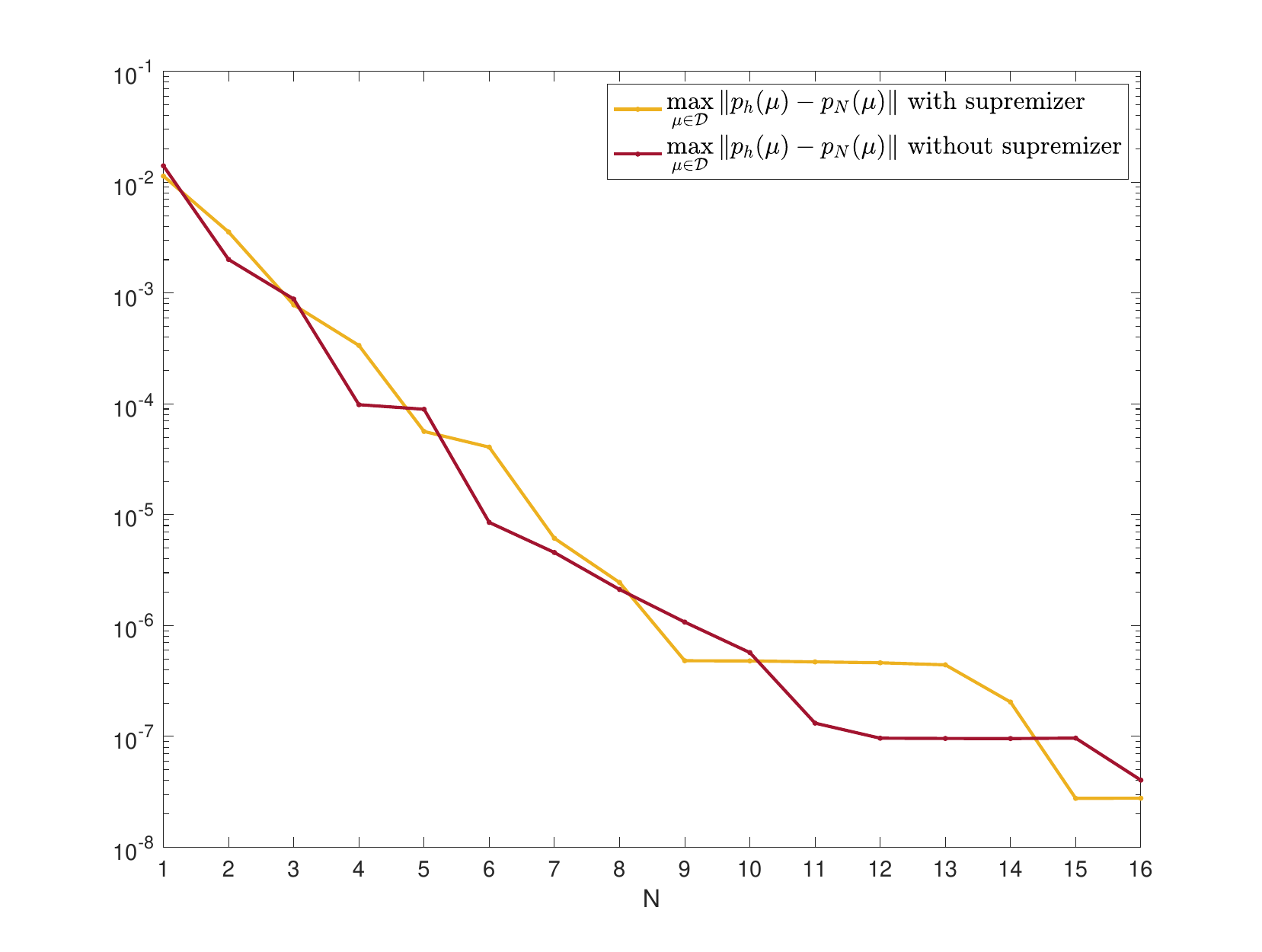}
\caption{Test 1. Velocity (left) and pressure (right) maximum error during the Greedy algorithm.}\label{chap:LPS::fig:error_vel_pres_TH}
\end{figure}

In Fig. \ref{chap:VMS::fig:DeltaCav_TH}, we show the comparison of the exact error for some snapshots and the \textit{a posteriori} error bound estimator considering either including the supremizers in the velocity spaces (left) or not (right). We can observe that the effectivity in this case is about one order of magnitude. The estimator is somewhat less smooth as function of the parameter $\mu$ for the LPS-VMS method.

\begin{figure}[h]
\centering
\includegraphics[width=0.45\linewidth]{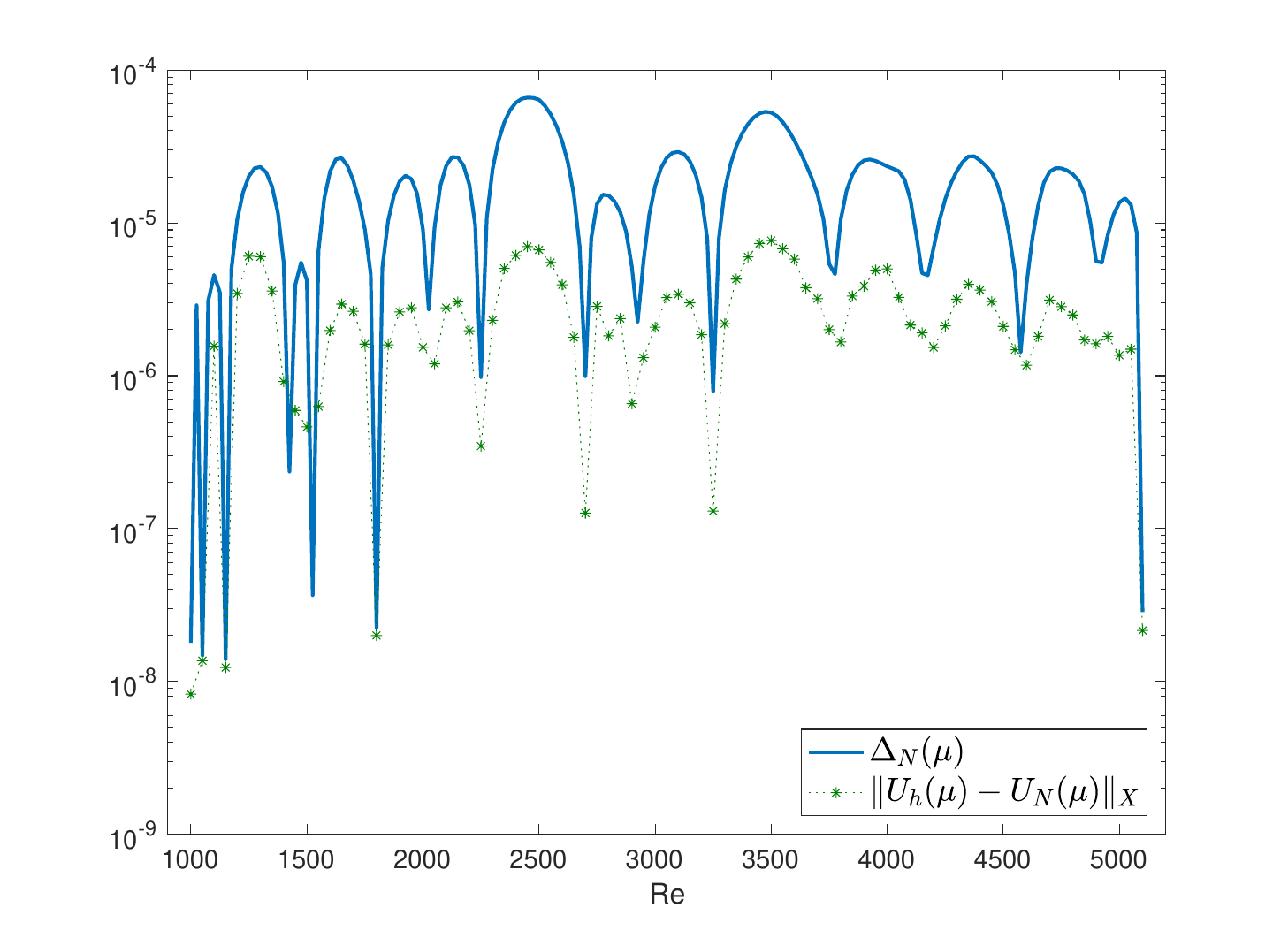}
\includegraphics[width=0.45\linewidth]{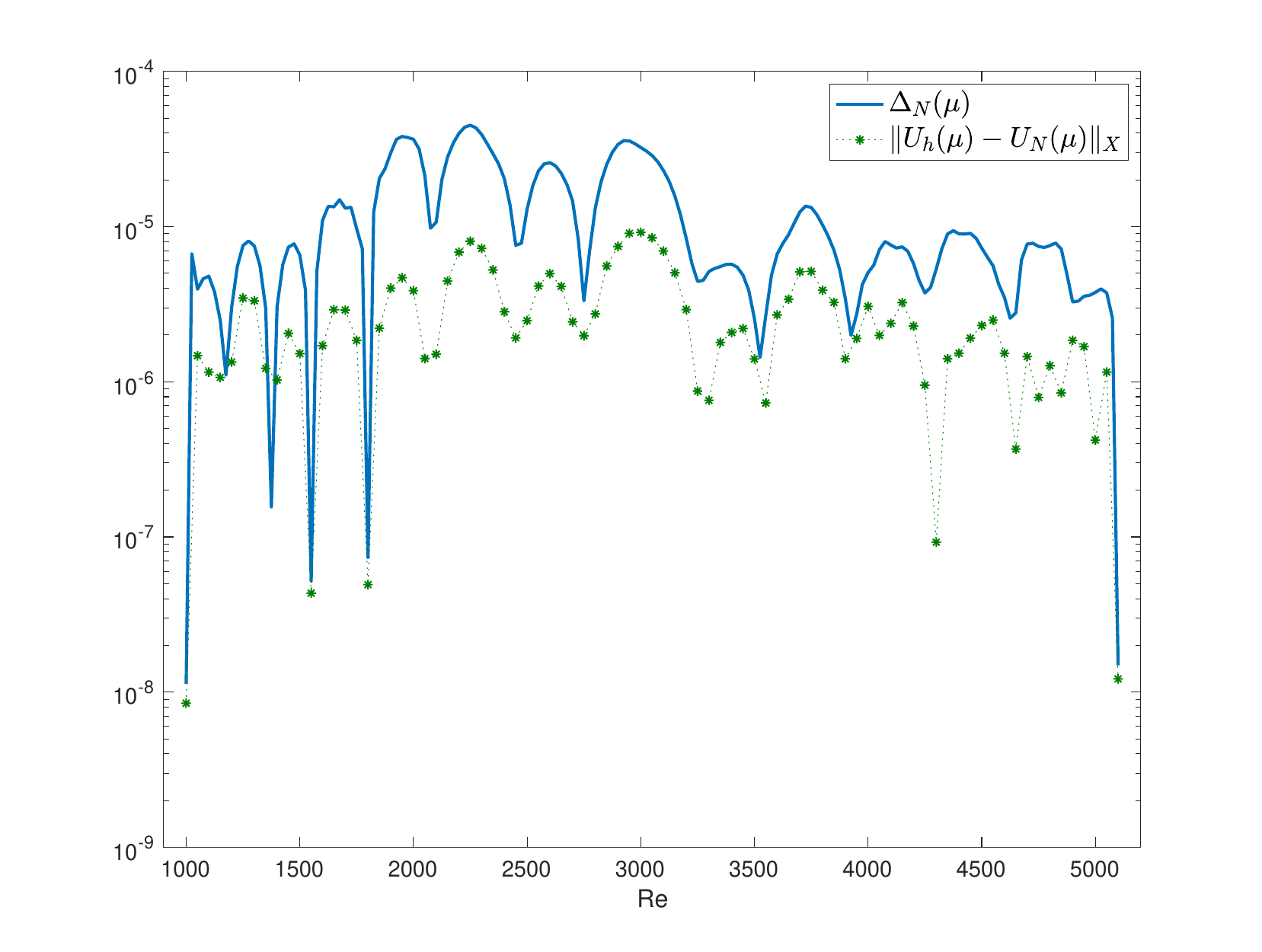}
\caption{Error and \textit{a posteriori} error bound estimator for $N=16$, with (left) and without pressure supremizers (right).}\label{chap:VMS::fig:DeltaCav_TH}
\end{figure}

Finally, in Table \ref{chap:VMS::tab:resCav_TH} we summarize results obtained for several values of \black{Reynolds} numbers, where we can observe that the absolute error between FE and RB solution is $10^{-6}$ for the velocity and $10^{-8}$ for pressure. Moreover, we can observe that the speedup of the solution is several \black{orders} of magnitude, with higher speedups (increase of nearly 50\%) and somewhat larger errors for the strategy of not including the pressure \textit{supremizers} in the reduced pressure spaces.

\begin{table}[h]
\begin{small}
\[\hspace{-0.1cm}
\begin{tabular}{l|cccc}
\hline
\multicolumn{5}{c}{\text{\textbf{Without supremizer}}}\\
\hline
Data &$\mu=1610$&$\mu=2751$&$\mu=3886$&$\mu=4521$\\
\hline
$T_{FE}$&1372.65s & 1513.53s & 3184.17s &3608.32s\\ 
$T_{online}$& 0.57s& 0.65s& 0.7s&0.77s\\
\hline
speedup& 2407 & 2327& 4548 &4685\\
\hline
$\|\uk_h-\uk_N\|_T$&$2.01\cdot10^{-6}$&$1.98\cdot10^{-6}$ & $2.04\cdot10^{-6}$ &$2.64\cdot10^{-6}$\\
$\|p_h-p_N\|_{0,2,\Omega}$&$1.51\cdot10^{-8}$&$2.03\cdot10^{-8}$ & $9.73\cdot10^{-9}$  &$2.05\cdot10^{-8}$\\
\hline
\end{tabular}\]
\[\hspace{-0.1cm}
\begin{tabular}{l|cccc}
\hline
\multicolumn{5}{c}{\text{\textbf{With supremizer}}}\\
\hline
Data &$\mu=1610$&$\mu=2751$&$\mu=3886$&$\mu=4521$\\
\hline
$T_{FE}$&1372.65s & 1513.53s & 3184.17s &3608.32s\\ 
$T_{online}$& 0.85s& 0.97s& 1.08s&1.11s\\
\hline
speedup& 1596 & 1544& 2947 &3224\\
\hline
$\|\uk_h-\uk_N\|_T$&$2.21\cdot10^{-6}$&$2.85\cdot10^{-6}$ & $3.83\cdot10^{-7}$ &$1.98\cdot10^{-6}$\\
$\|p_h-p_N\|_{0,2,\Omega}$&$1.09\cdot10^{-8}$&$1.37\cdot10^{-8}$ & $7.08\cdot10^{-9}$  &$1.10\cdot10^{-8}$\\
\hline
\end{tabular}\]
\caption{Computational time for FE solution and RB online phase, with the speedup and the error, \black{for $\mathbb{P}_2-\mathbb{P}_1$.}}
\label{chap:VMS::tab:resCav_TH}
\end{small}
\end{table}

\subsection{Test 2: $\mathbb{P}_2-\mathbb{P}_2$ velocity-pressure finite element}

For this test, we need $M_1=52$ basis for the eddy viscosity approximation and $M_2=28$ basis for the pressure stabilisation constant approximation to reach a \black{prescribed} tolerance of $\veps_{EIM}=5\cdot10^{-5}$ for both terms. In Fig. \ref{chap:VMS::fig:EIMCav} we show the evolution of the error for both non-linear terms.

\begin{figure}[h]
\centering
\includegraphics[width=0.45\linewidth]{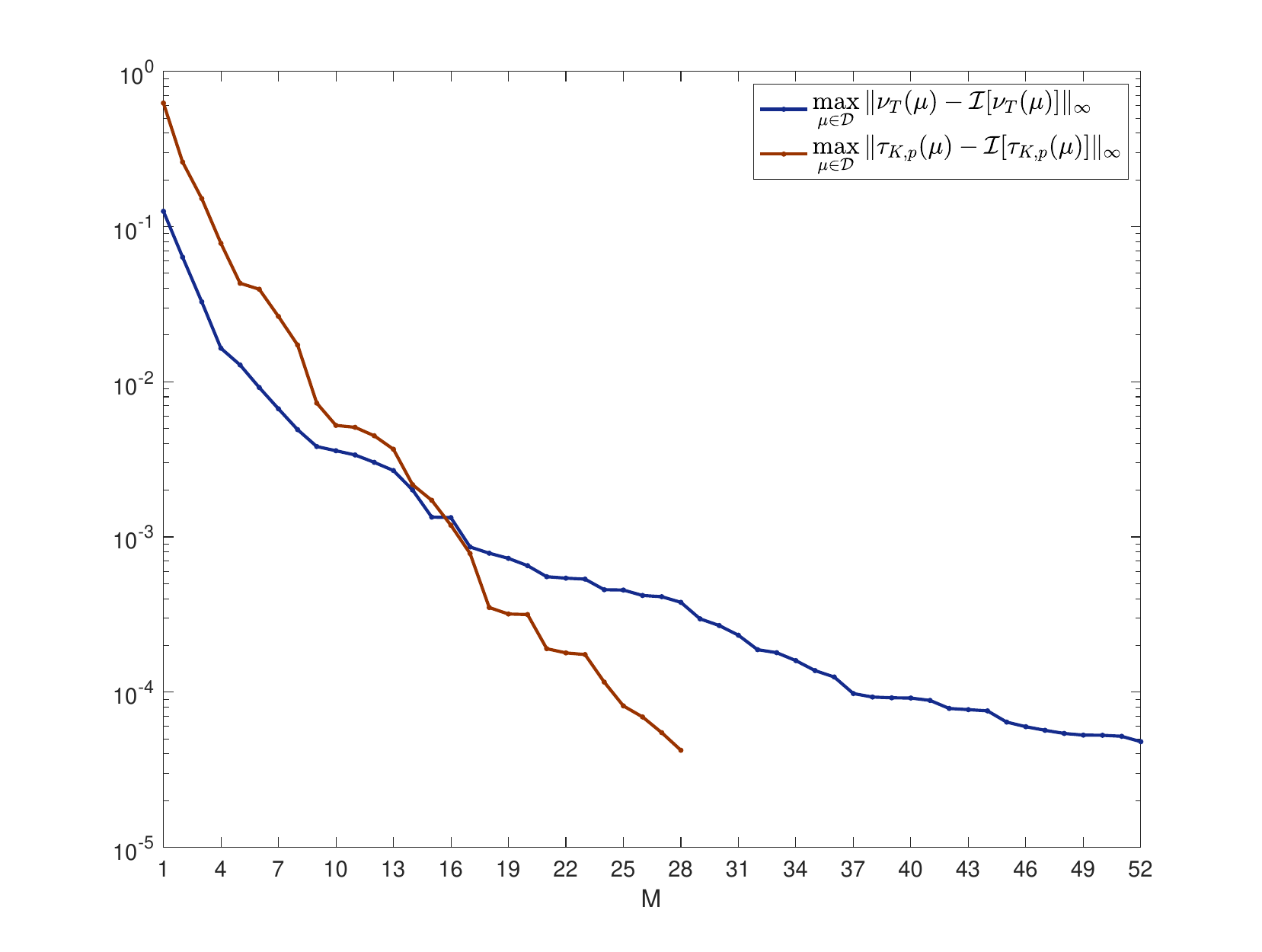}
\caption{Evolution of the error in the EIM.}\label{chap:VMS::fig:EIMCav}
\end{figure}

For the Greedy algorithm to determine the reduced velocity and pressure spaces, in this case we prescribe the same tolerance of $\veps_{RB}=7\cdot10^{-5}$. This tolerance is reached for $N=16$ basis functions when we do not consider supremizers, while in the case of considering the supremizer, the tolerance is reached for $N=17$. In Fig. \ref{chap:VMS::fig:EstErrCav} we show the evolution of the \textit{a posteriori} error bound estimator during the Greedy algorithm for both cases. We can observe, again, that the evolution of the maximum value for the error estimator either considering supremizers (left) or not (right) are equivalent. Although we do not know wether the reduced problem without su\-pre\-mi\-zers is stable, all computations take place as if it was.

\begin{figure}[h]
\centering
\includegraphics[width=0.45\linewidth]{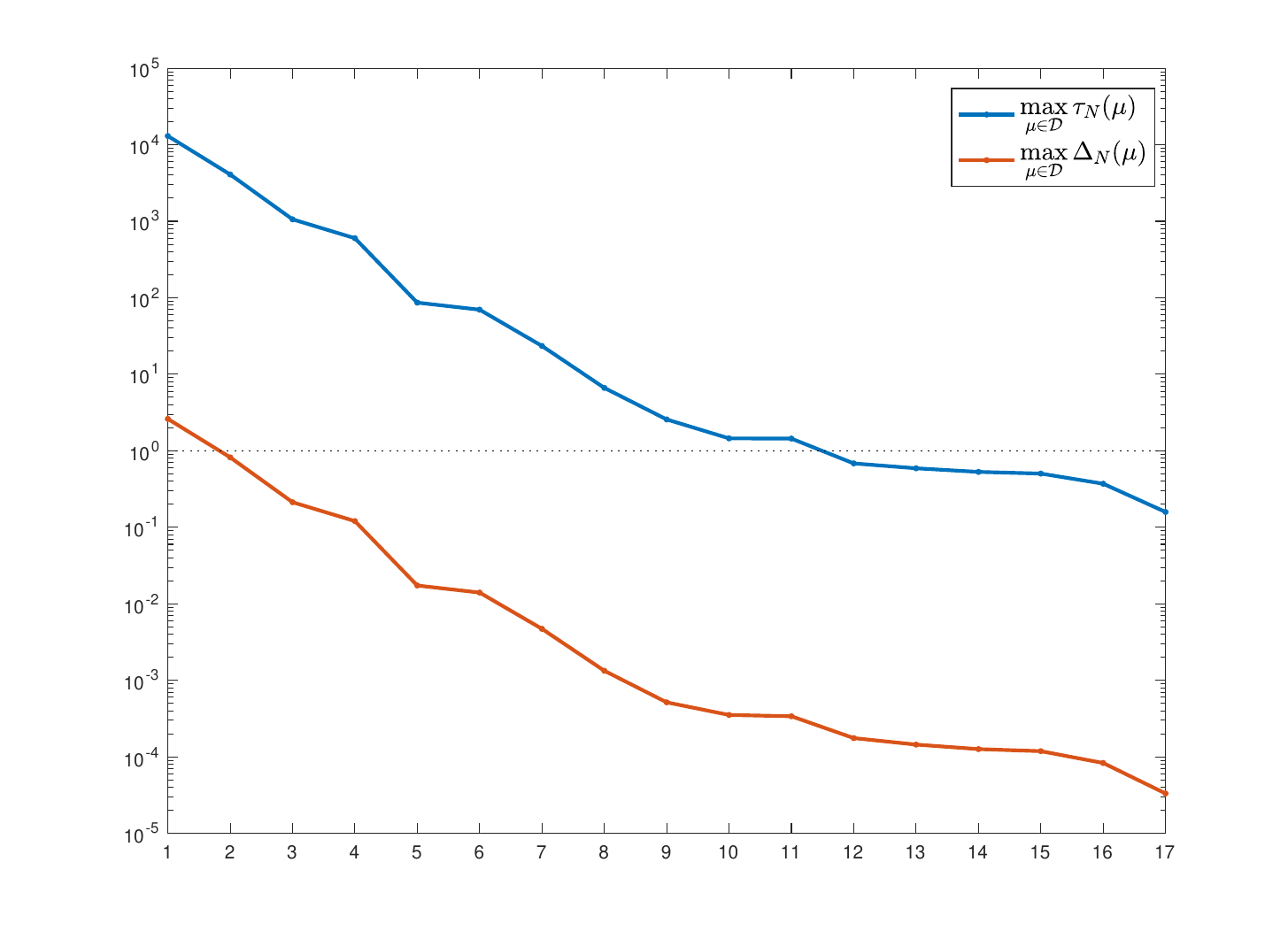}
\includegraphics[width=0.45\linewidth]{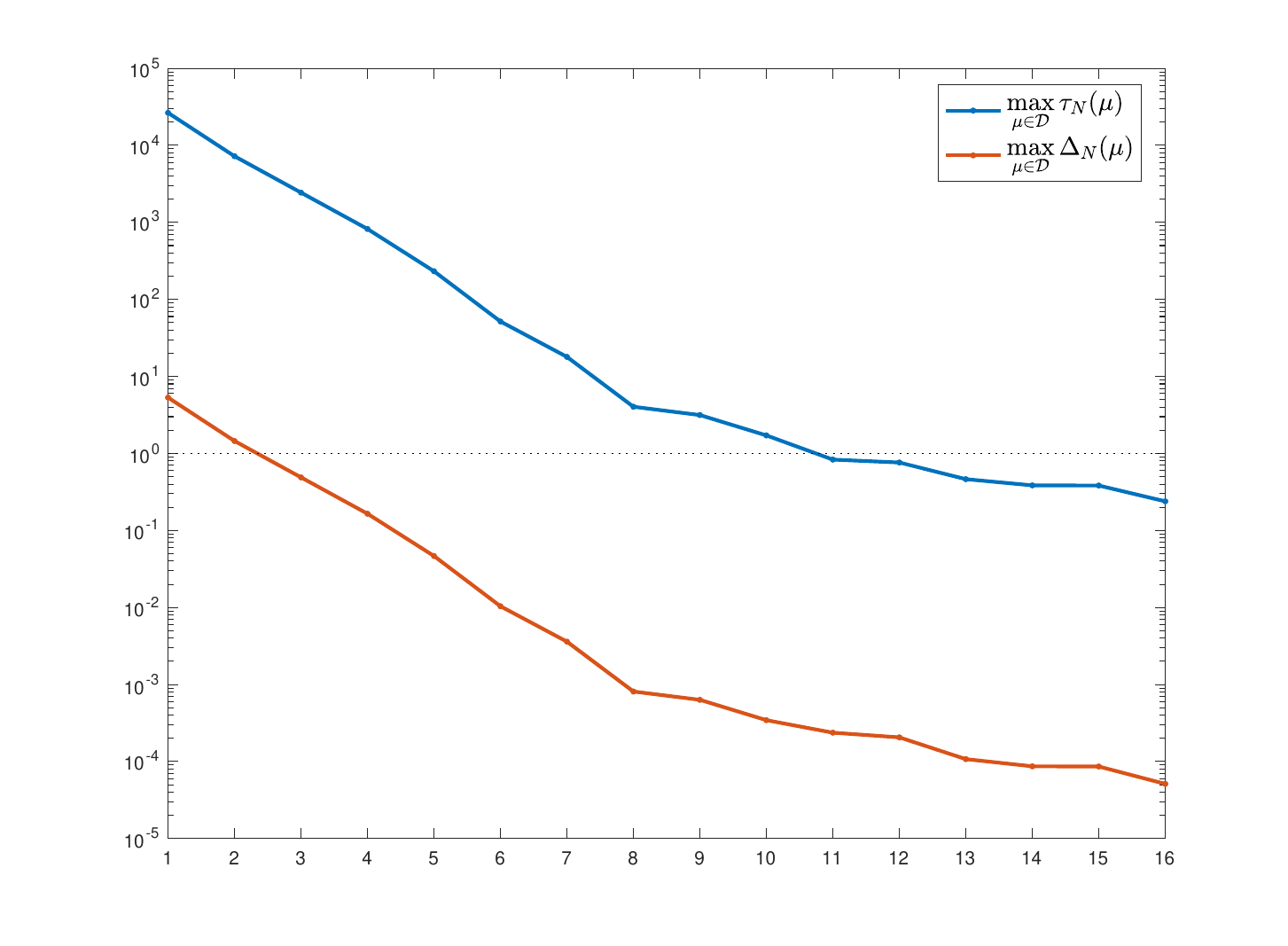}
\caption{Test 2. Evolution of the error in the Greedy algorithm with (left) and without pressure supremizers (right).}\label{chap:VMS::fig:EstErrCav}
\end{figure}

We can observe in Fig. \ref{chap:LPS::fig:error_vel_pres} that for both the velocity and pressure the maximum error value in each iteration of the Greedy algorithm, for both velocity and pressure, are quite close, as in the Taylor-Hood finite element test.

\begin{figure}[h]
\centering
\includegraphics[width=0.45\linewidth]{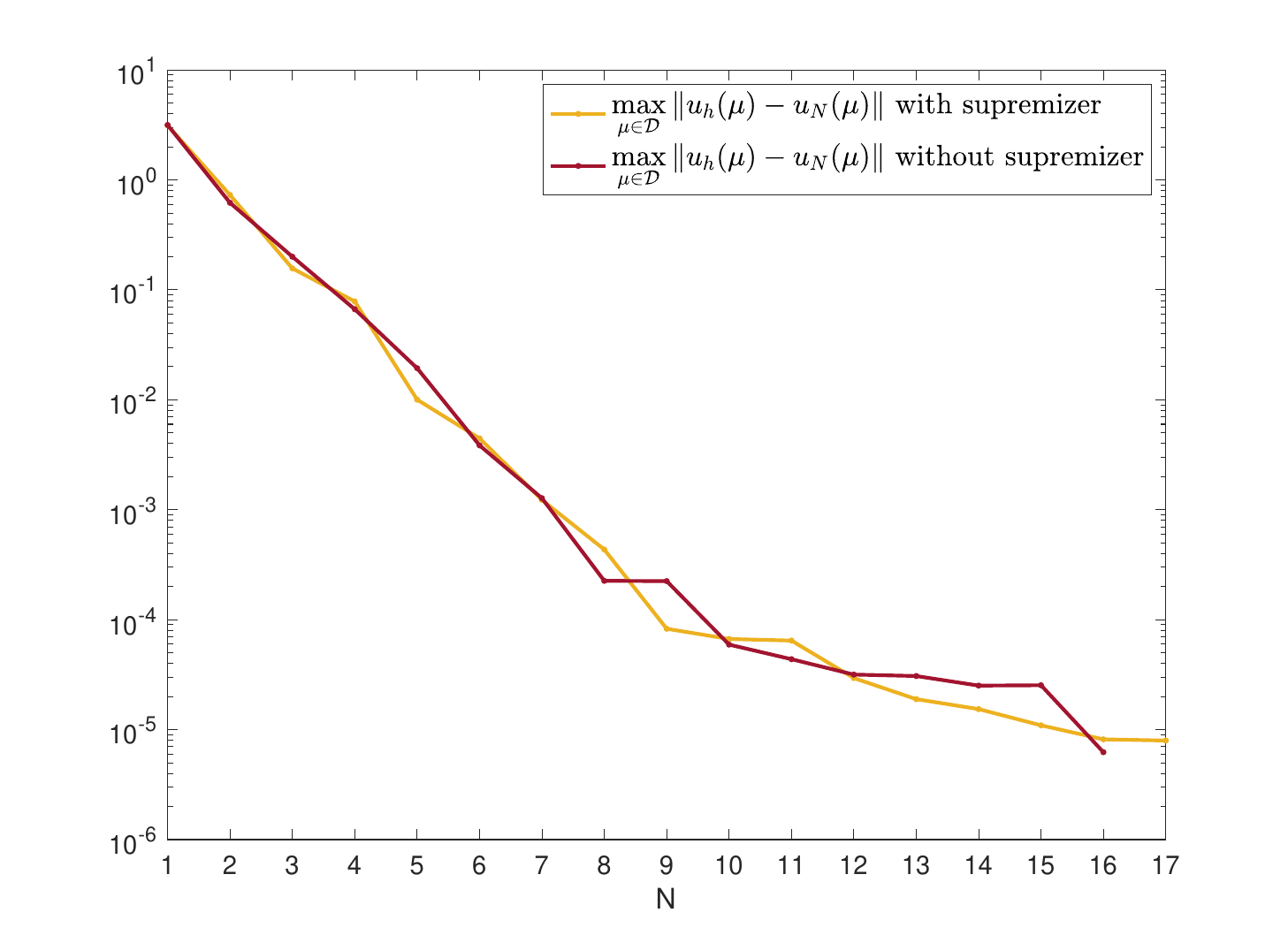}
\includegraphics[width=0.45\linewidth]{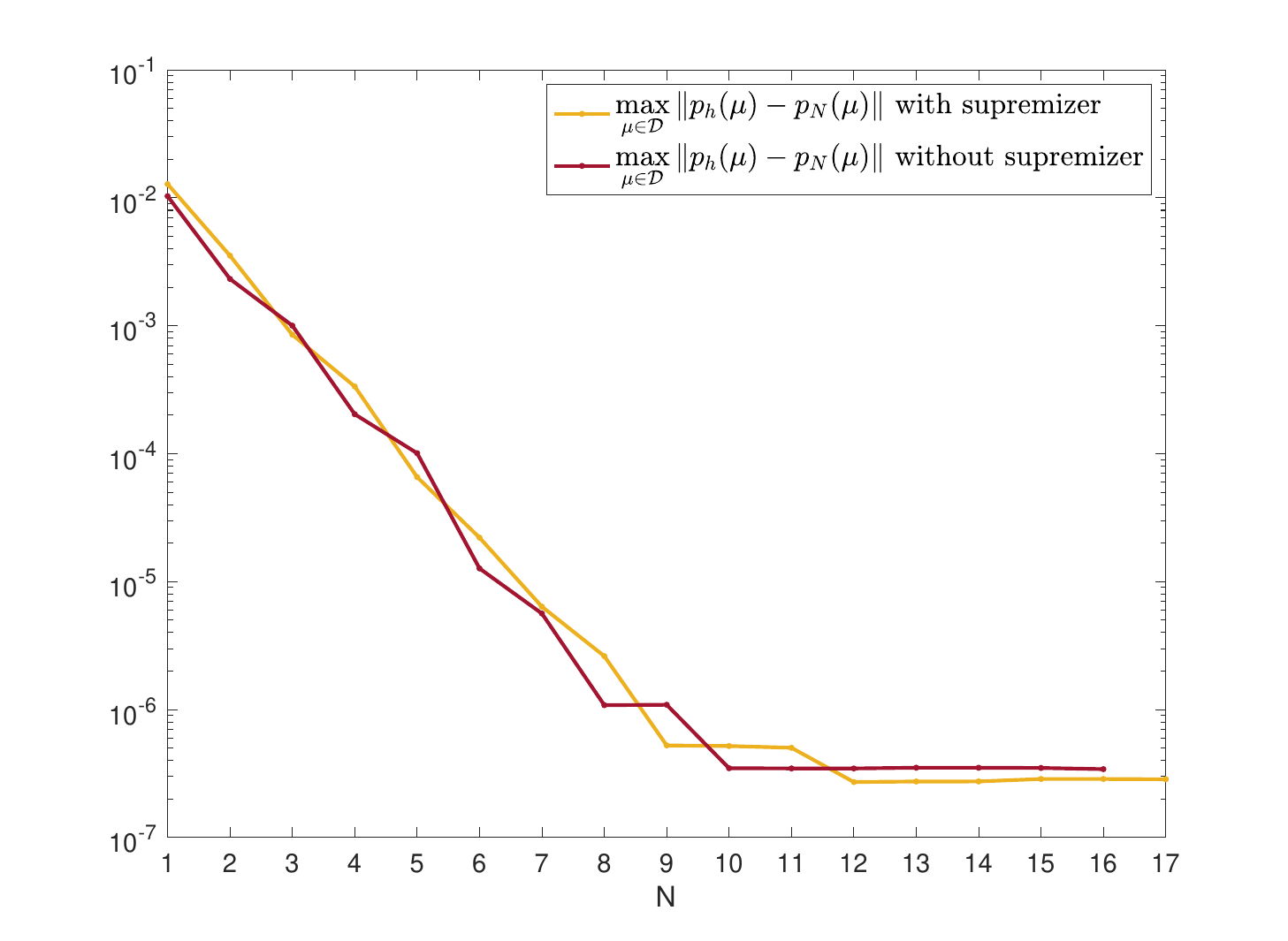}
\caption{Test 2. Velocity (left) and pressure (right) maximum error during the Greedy algorithm.}\label{chap:LPS::fig:error_vel_pres}
\end{figure}

In Fig. \ref{chap:VMS::fig:DeltaCav} we show the comparison of the exact error for some snapshots and the \textit{a posteriori} error bound estimator considering either including the supremizers in the velocity spaces (left) or not (right).

\begin{figure}[h]
\centering
\includegraphics[width=0.45\linewidth]{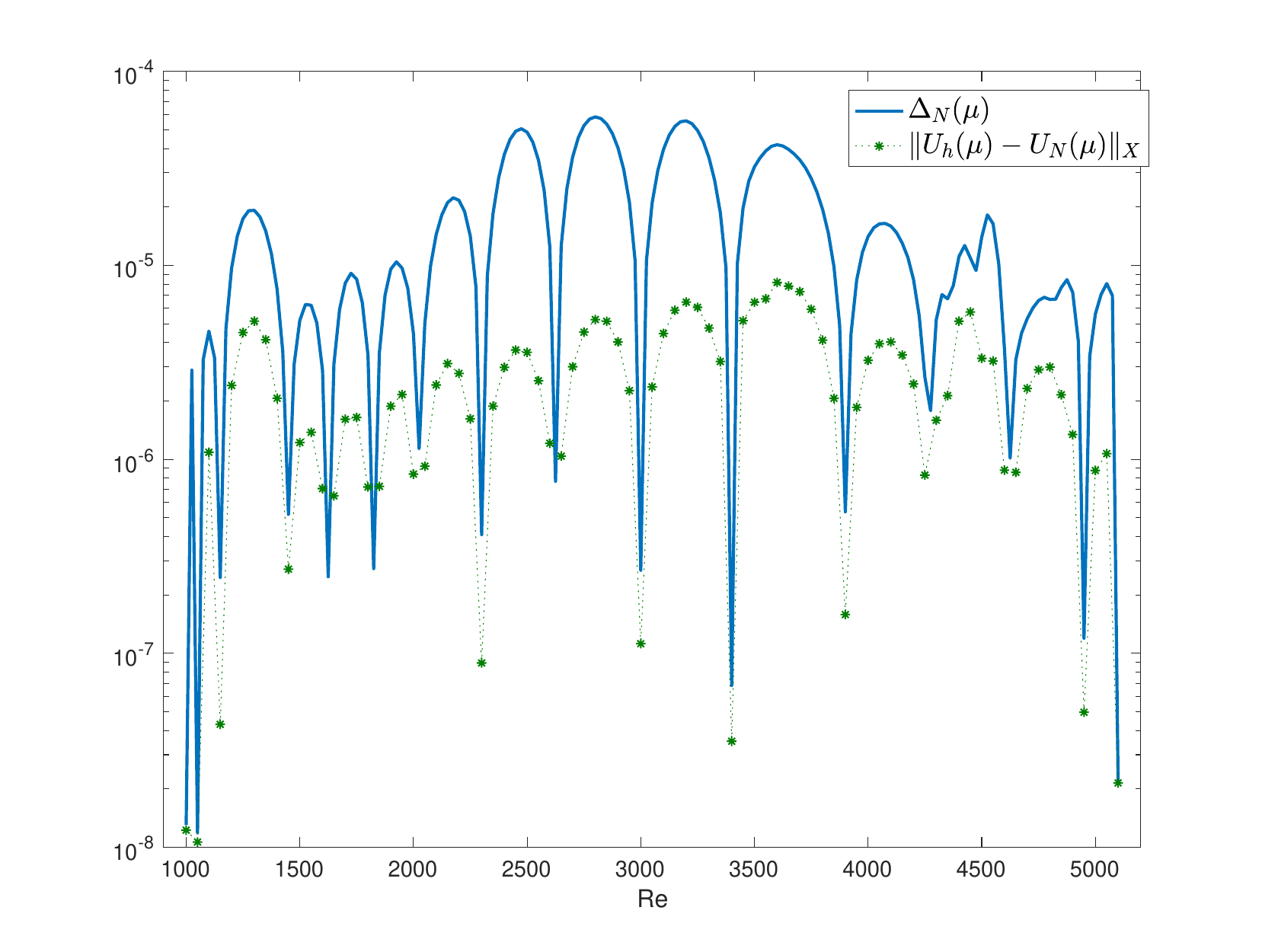}
\includegraphics[width=0.45\linewidth]{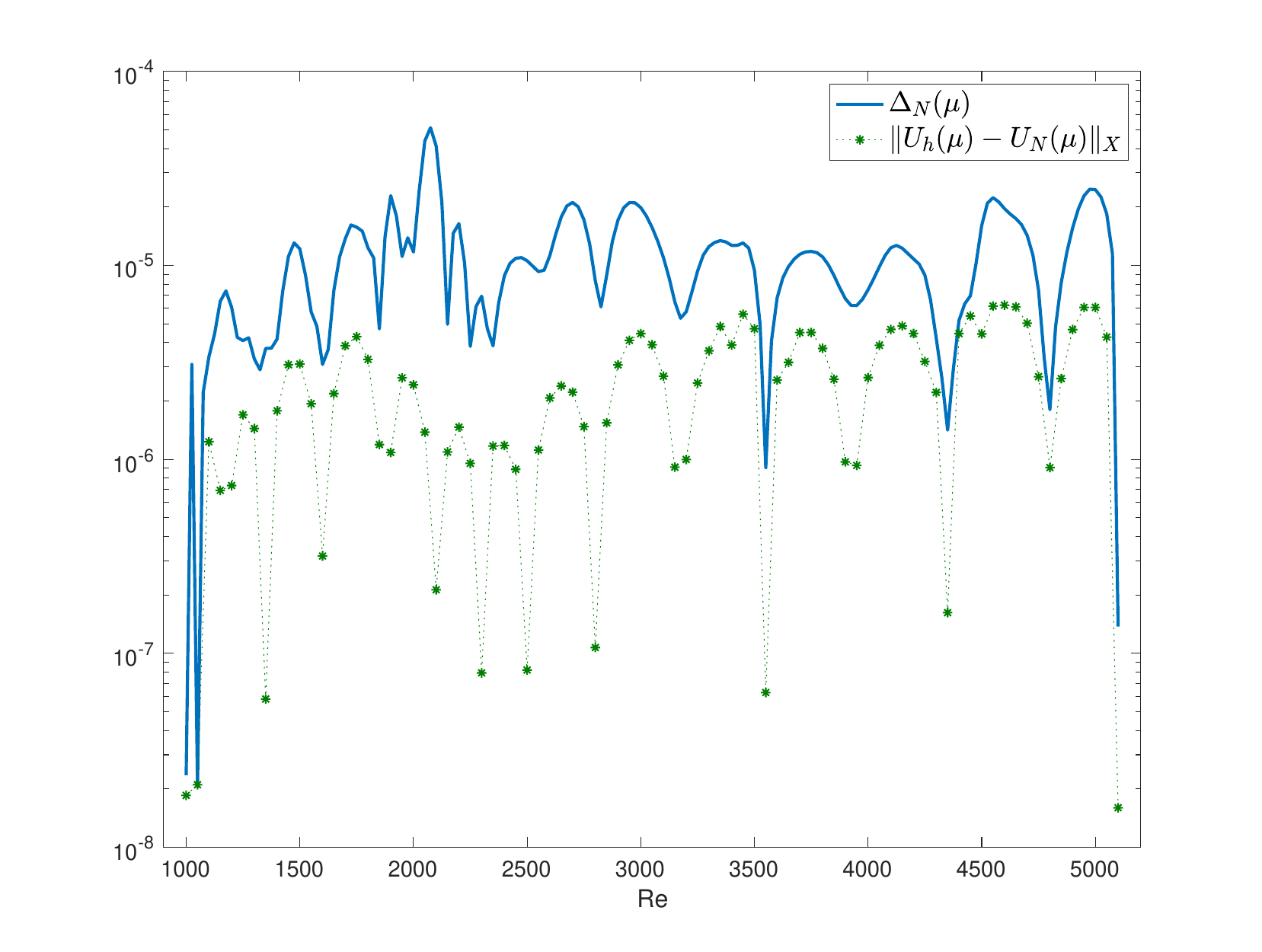}
\caption{Test 2. Error and \textit{a posteriori} error bound estimator for $N=16$, with (left) and without pressure supremizers (right).}\label{chap:VMS::fig:DeltaCav}
\end{figure}

Finally, in Table \ref{chap:VMS::tab:resCav} we summarize results obtained for several values of \black{Reynolds} numbers, where we can observe that the absolute error between FE and RB solution is around order of $10^{-6}$ for the velocity and $10^{-8}$ for pressure. Moreover, we can observe that the speedup of the solution is of several orders of magnitude, with increase of nearly 50\%, and somewhat larger errors, when no supremizers are used.

\begin{table}[h]
\begin{small}
\[\hspace{-0.1cm}
\begin{tabular}{l|cccc}
\hline
\multicolumn{5}{c}{\text{\textbf{Without supremizer}}}\\
\hline
Data &$\mu=1610$&$\mu=2751$&$\mu=3886$&$\mu=4521$\\
\hline
$T_{FE}$&2259.04s & 3008.81s & 5574.5s &6171.41s\\ 
$T_{online}$& 0.70s& 0.75s& 0.78s&0.81s\\
\hline
speedup& 3227 & 4011& 7146 &7619\\
\hline
$\|\uk_h-\uk_N\|_T$&$5.4\cdot10^{-7}$&$1.44\cdot10^{-6}$ & $1.49\cdot10^{-6}$ &$5.44\cdot10^{-6}$\\
$\|p_h-p_N\|_{0,2,\Omega}$&$1.34\cdot10^{-8}$&$3.35\cdot10^{-8}$ & $1.55\cdot10^{-7}$  &$2.25\cdot10^{-8}$\\
\hline
\end{tabular}\]
\[\hspace{-0.1cm}
\begin{tabular}{l|cccc}
\hline
\multicolumn{5}{c}{\text{\textbf{With supremizer}}}\\
\hline
Data &$\mu=1610$&$\mu=2751$&$\mu=3886$&$\mu=4521$\\
\hline
$T_{FE}$&2259.04s & 3008.81s & 5574.5s &6171.41s\\ 
$T_{online}$& 0.99s& 1.03s& 1.26s&1.23s\\
\hline
speedup& 2267 & 2895& 4391 &5016\\
\hline
$\|\uk_h-\uk_N\|_T$&$4.94\cdot10^{-7}$&$4.51\cdot10^{-6}$ & $6.52\cdot10^{-7}$ &$3.52\cdot10^{-6}$\\
$\|p_h-p_N\|_{0,2,\Omega}$&$5.23\cdot10^{-8}$&$3.22\cdot10^{-8}$ & $6.38\cdot10^{-8}$  &$8.45\cdot10^{-8}$\\
\hline
\end{tabular}\]
\caption{Computational time for FE solution and RB online phase, with the speedup and the error, for $\mathbb{P}_2-\mathbb{P}_2$.}
\label{chap:VMS::tab:resCav}
\end{small}
\end{table}

\section{Conclusions}
In this paper we have developed a VMS-Smagorinsky reduced basis model with local projection stabilisation (LPS) to stabilise the pressure discretization. We have proved the stability of the LPS stabilised reduced problem for piecewise affine pressure discretisations. We \black{have used} the EIM to approximate the eddy viscosity and the pressure stabilisation coefficient in order to reduce the non-linearities that appear in both terms, letting us to efficiently store the tensors associated to both the VMS-Smagorinsky term and the LPS pressure term. Also, we have developed an \textit{a posteriori} error estimator for this problem via the BRR theory. 

Finally, we have presented some numerical results of the RB method cons\-truc\-ted by means of a Greedy algorithm, based upon this error estimator, applied to the 2D lid-driven cavity problem.  We \black{have used} both Taylor-Hood and $\mathbb{P}_2-\mathbb{P}_2$ pairs of velocity-pressure discretizations of the full order model. Our theory proves that the Taylor-Hood discretisation reduced problem including LPS stabilisation of the pressure is stable. We have observed that in both strategies (including or not pressure supremizers in the reduced velocity spaces), the errors obtained for the velocity and pressure reduced solutions are quite close. Further, the computational speedup is higher, with increases of nearly 30\%, when we do not consider the enrichment with pressure supremizers. We thus have built a RB method with enhanced reduction and speedup, while keeping error levels close to using pressure supremizers.
\section*{Acknowledgements}
This work has been partially supported by the Andalusian Government - Feder Fund Project US-1254587, that provided the funding to hire Dr. Delgado for some time, and the Spanish Government - Feder Fund Project RTI2018-093521-B-C31.
\bibliographystyle{siam}
\bibliography{Ref_RBM.bib}

\end{document}